\documentclass{amsart}
\usepackage{amssymb,amsfonts}
\usepackage{enumerate}
\usepackage{enumitem}
\usepackage{mathrsfs}
\usepackage{graphicx}
\usepackage{mathtools}
\usepackage{leftidx}
\usepackage{tikz-cd}
\usepackage{amsmath,calligra,mathrsfs}
\usepackage{bm}
\usepackage{url}
\newtheorem{thm}{Theorem}[section]
\newtheorem{cor}[thm]{Corollary}
\newtheorem{prop}[thm]{Proposition}
\newtheorem{lem}[thm]{Lemma}

\theoremstyle{definition}
\newtheorem{defn}[thm]{Definition}

\newtheorem{exmp}[thm]{Example}

\newtheorem{conv}[thm]{Convention}

\theoremstyle{remark}
\newtheorem{rem}[thm]{Remark}

\makeatletter
\let\c@equation\c@thm
\makeatother
\numberwithin{equation}{section}

\def\bthm{\begin{thm}}
\def\ethm{\end{thm}}
\def\blm{\begin{lem}}
\def\elm{\end{lem}}
\def\bdf{\begin{defn}}
\def\edf{\end{defn}}
\def\bpf{\begin{proof}}
\def\epf{\end{proof}}
\def\bpp{\begin{prop}}
\def\epp{\end{prop}}
\def\bcor{\begin{cor}}
\def\ecor{\end{cor}}
\def\brm{\begin{rem}}
\def\erm{\end{rem}}
\def\beg{\begin{exmp}}
\def\eeg{\end{exmp}}
\def\bD{\mathbb{D}}
\def\bE{\mathbb{E}}

\def\bG{\mathbb{G}}
\def\bH{\mathbb{H}}

\def\bN{\mathbb{N}}

\def\bQ{\mathbb{Q}}

\def\bZ{\mathbb{Z}}
\def\cA{\mathcal{A}}
\def\cB{\mathcal{B}}
\def\cC{\mathcal{C}}
\def\cD{\mathcal{D}}
\def\cE{\mathcal{E}}
\def\cF{\mathcal{F}}
\def\cG{\mathcal{G}}
\def\cH{\mathcal{H}}

\def\cK{\mathcal{K}}
\def\cL{\mathcal{L}}
\def\cM{\mathcal{M}}
\def\cN{\mathcal{N}}
\def\cO{\mathcal{O}}
\def\cP{\mathcal{P}}

\def\cR{\mathcal{R}}

\def\scH{\mathscr{H}}
\def\scI{\mathscr{I}}
\def\scJ{\mathscr{J}}
\def\scK{\mathscr{K}}
\def\scL{\mathscr{L}}

\def\scO{\mathscr{O}}

\newcommand{\raq}{\,\rightarrow \,}
\newcommand{\laq}{\,\leftarrow \,}

\newcommand{\xraq}[2][]{\, \xrightarrow[#1]{#2} \,}
\newcommand{\xlaq}[2][]{\, \xleftarrow[#1]{#2} \,}

\newcommand{\ra}{\rightarrow}
\newcommand{\la}{\leftarrow}
\newcommand{\rinto}{\hookrightarrow}

\newcommand{\xra}[2][]{\xrightarrow[#1]{#2}}
\newcommand{\xla}[2][]{\xleftarrow[#1]{#2}}
\newcommand{\xrinto}[2][]{\xhookrightarrow[#1]{#2}}

\newcommand{\ie}{{\it i.e.}}
\newcommand{\eg}{{\it e.g.}}

\newcommand{\Mod}{{\rm Mod}}

\newcommand{\Ch}{{\rm Ch}}
\newcommand{\Tr}{{\rm Tr}}

\newcommand{\Tor}{{\rm Tor}}

\newcommand{\EExt}{\bE{\rm xt}}

\newcommand{\op}{{\rm op}}

\newcommand{\id}{{\rm id}}

\newcommand{\Hom}{{\rm Hom}}
\newcommand{\Homcom}{\underline{{\rm Hom}}}

\newcommand{\cHom}{\mathscr{H}\text{\kern -3pt {\calligra\large om}}\,}

\newcommand{\RHom}{{\bm R}{\rm Hom}}
\newcommand{\RHomcom}{{\bm R} \underline{{\rm Hom}}}

\newcommand{\RcHom}{{\bm R}\cHom}
\newcommand{\RcHomcom}{{\bm R}\underline{\cHom}}

\newcommand{\Spec}{{\rm Spec}}

\newcommand{\Proj}{{\rm Proj}}
\newcommand{\PProj}{\mathpzc{Proj}}

\newcommand{\WDiv}{{\rm WDiv}}

\newcommand{\Ex}{{\rm Ex}}

\newcommand{\QCoh}{{\rm QCoh}}
\newcommand{\Coh}{{\rm Coh}}

\newcommand{\Dsuit}{\cD^{\tiny \mbox{$\spadesuit $}}}

\newcommand{\Dbcoh}{\cD^b_{{\rm coh}}}

\newcommand{\Dqcoh}{\cD_{{\rm qcoh}}}

\newcommand{\Dmcoh}{\cD^{-}_{{\rm coh}}}

\newcommand{\Dperf}{\cD_{{\rm perf}}}

\newcommand{\coh}{{\rm coh}}

\newcommand{\Pic}{{\rm Pic}}

\newcommand{\Gr}{{\rm Gr}}
\newcommand{\GrA}{{\rm Gr}(A)}
\newcommand{\gr}{{\rm gr}}

\newcommand{\QpGr}{{\rm Q}^+{\rm Gr}}

\newcommand{\qpgr}{{\rm q}^+{\rm gr}}

\newcommand{\QpdGr}{{\rm Q}^+_{(d)}{\rm Gr}}

\newcommand{\Torp}{{\rm Tor}^+}

\newcommand{\IIT}{I^{\infty}\text{-}{\rm Tor}}

\newcommand{\Icomp}{I\text{-}{\rm comp}}

\newcommand{\RGam}{{\bm R}\Gamma}

\newcommand{\RGI}{{\bm R}\Gamma_{I}}

\newcommand{\Ce}{\check{\cC}}
\newcommand{\CI}{\check{\cC}_I}

\newcommand{\LLI}{{\bm L}\Lambda_{I}}

\newcommand{\ITR}{I\text{-}{\rm triv}}
\newcommand{\IpTR}{I^+\text{-}{\rm triv}}
\newcommand{\ImTR}{I^-\text{-}{\rm triv}}

\newcommand{\cIpTR}{\scI^+\text{-}{\rm triv}}

\newcommand{\lc}{{\rm lc}}

\DeclareMathAlphabet{\mathpzc}{OT1}{pzc}{m}{it}

\title{Homological flips and homological flops}
\author{Wai-Kit Yeung}
\address{Department of Mathematics,
	Indiana University, Bloomington, IN 47405, USA}
\email{yeungw@iu.edu}
\begin{document}

\begin{abstract}
We introduce a notion of homological flips and homological flops. The former includes the class of all flips between Gorenstein normal varieties; while the latter includes the class of all flops between Cohen-Macaulay normal varieties whose contracted variety is quasi-Gorenstein. 
Our main theorem shows that certain local cohomology complexes are dual to each other under homological flips/flops. 
We give some preliminary applications of this duality to relate the derived categories under flip/flop. Further applications are in \cite{Yeu20b}.
\end{abstract}

\maketitle


\tableofcontents

\section{Introduction}

A large amount of information about the geometry of a variety is encoded in its derived category. This rich interplay between the homological algebra of derived categories and the geometry of algebraic varieties is especially apparent in birational geometry. 
Indeed, in the seminal paper \cite{BO95}, Bondal and Orlov put forth the following principle:
\begin{equation}  \label{CMMP_principle}
\parbox{40em}{The minimal model program (MMP) amounts to minimizing the derived category in a given birational class.}
\end{equation}

Bondal and Orlov have provided some evidence for this principle in \cite{BO95}. Since then, more and more results of the same spirit has been established by the work of many other mathematicians
(see, e.g., the surveys \cite{Kaw09, Kaw18} and the references therein).
A closely related phenomenon is also given by the ``DK-conjecture" of Kawamata (see, e.g., \cite[Conjecture 1.2]{Kaw18}), which could be viewed as a variation of the same theme.

Recall that there are three basic birational operations in the minimal model program, namely divisorial contractions, flips, and flops. Thus, principle \eqref{CMMP_principle} asserts in particular that
\begin{equation}  \label{CMMP_principle_2}
\parbox{40em}{The derived category should shrink under divisorial contractions and flips, and it should remain equivalent under flops. }
\end{equation}

There is, however, an ambiguity of what it means by the ``derived category", and what it means for the derived category to ``shrink". For example, the derived category of a variety $X$ may refer to either $\Dbcoh(X)$ or $\Dperf(X)$, which coincide if and only if $X$ is regular; and shrinking may refer to a fully faithful functor, a localization, or a semi-orthogonal decomposition. In this introduction, we shall take the liberty to switch between these alternative viewpoints.

If we interpret the ``derived category" to mean $\Dperf(X)$, and ``shrinking'' to mean fully faithful embedding
then it is easy to see that the derived category often shrinks under divisorial contractions.
Indeed, varieties in the MMP often have rational singularities, and hence any birational contraction $f : X \ra Y$ between them satisfies 
$f_*(\cO_X) = \cO_Y$ and ${\bm R}^i f_*(\cO_X) = 0$ for all $i > 0$, which is equivalent to the functor  ${\bm L}f^* : \Dperf(Y) \ra \Dperf(X)$ being fully faithful.

For flips and flops, such a direct comparison between the derived categories is absent because there is no direct morphism between the varieties under flips and flops. The main purpose of this paper is to study the change in the derived categories under flips and flops.


Consider a log flip between normal varieties (see \eqref{Xm_Y_Xp_2} below for definition)
\begin{equation}  \label{log_flip_intro}
\begin{tikzcd} [row sep = 0]
X^-  \ar[rd,"\pi^-"] & &  X^+ \ar[ld, "\pi^+"'] \\
& Y
\end{tikzcd}
\end{equation}
By definition of a log flip, 
there are Weil divisors $D^-$ on $X^-$ and $D^+$ on $X^+$, strict transforms of each other, such that
\begin{equation}  \label{divisor_D_pm_intro}
\parbox{40em}{(1) $-D^-$ is $\bQ$-Cartier and $\pi^-$-ample; \\
(2) $D^+$ is $\bQ$-Cartier and $\pi^+$-ample.}
\end{equation}
If we denote by $D_Y$ their common strict transform to $Y$, and let $\cA$ be the sheaf of $\bZ$-graded algebras $\cA = \bigoplus_{i \in \bZ} \cO_Y(iD_Y)$, then we have $X^- = \Proj_Y(\cA_{\leq 0})$ and $X^+ = \Proj_Y(\cA_{\geq 0})$.
For simplicity of language, we will assume for now that $Y$ is affine. Thus, the log flip \eqref{log_flip_intro} is completely determined by a finitely generated $\bZ$-graded algebra $A$.

Readers familiar with variations of GIT quotients will notice that this is an instance of a wall-crossing. Indeed, let $W = \Spec \, A$, with a $\bG_m$-action coming from the $\bZ$-grading, and let $L^-$ (resp. $L^+$) be the trivial line bundle with $\bG_m$-linearizations given by the $-1$ (resp. $+1$) character of $\bG_m$. Then the semistable loci are given by $W^{ss}(L^{\pm}) = W \setminus V(I^{\pm})$, where $V(I^{\pm})$ are the $\bG_m$-invariant closed subset defined by the graded ideals $I^+ := A_{>0} \cdot A$ and $I^- := A_{<0} \cdot A$ respectively.
The stack quotient by $\bG_m$ of the semistable loci are then the stacky projective spaces $\mathfrak{X}^{\pm} = \PProj^{\pm}(A) = [W^{ss}(L^{\pm}) / \bG_m]$, which are Deligne-Mumford stacks whose coarse moduli space are the varieties $X^{\pm}$.

In recent years, many results along the lines of \eqref{CMMP_principle} are established by techniques that relate the derived categories of the stacks $\mathfrak{X}^{\pm}$ under the wall-crossing (see, \eg, \cite{HL15, BFK19}). In the present context, this problem can be expressed in terms of the commutative algebra of the graded ring $A$.
%
For instance, let $f_1,\ldots,f_r$ be homogeneous elements in $A$ of positive degrees that generate $I^+$. 
For any graded module $M \in \Gr(A)$, 
define the \emph{\v{C}ech complex} in the usual way:
\begin{equation}  \label{Ce_f_intro}
\Ce_{I^+}(M)\, := \, \bigl[ \, \prod_{1 \leq i_0 \leq r} M_{f_{i_0}}  \xraq{-d^1}  
\prod_{1 \leq i_0 < i_1 \leq r} M_{f_{i_0}f_{i_1}}  \xraq{-d^2} \ldots \xraq{-d^{r-1}} 
M_{f_1\ldots f_r} \, \bigr]
\end{equation}

Denote by $j : \mathfrak{X}^+ \subset \mathfrak{X}$ the inclusion, then $\Ce_{I^+}(M)$ is the object in the derived category $\cD(\GrA) \simeq \cD(\QCoh(\mathfrak{X}))$ corresponding to ${\bm R}j_* j^*(M)$. In particular, the \v{C}ech complex sits in an exact triangle in the derived category $\cD(\GrA)$:
\begin{equation} \label{RGI_CI_seq_intro}
\ldots \raq \RGam_{I^+}(M) \xraq{\epsilon} M \xraq{\eta} \Ce_{I^+}(M) \xraq{\delta} \RGam_{I^+}(M)[1] \raq \ldots 
\end{equation}
where $\RGam_{I^+}(M) \in \cD(\GrA)$ is the local cohomology complex with respect to the graded ideal $I^+$. In fact, \eqref{RGI_CI_seq_intro} is the decomposition sequence associated to a semi-orthogonal decomposition
\begin{equation}  \label{torsion_SOD_intro}
\cD(\GrA) \, = \, \langle \, \cD_{\IpTR}(\GrA) \, , \,  \cD_{\Torp}(\GrA) \, \rangle
\end{equation}
where $\cD_{\IpTR}(\GrA)$ is equivalent to $\cD(\QCoh(\mathfrak{X}^+))$ via the functor ${\bm R}j_*$, and $\cD_{\Torp}(\GrA)$ consists of objects in  $\cD(\QCoh(\mathfrak{X}))$ supported along the unstable locus $V(I^+)$.

The \v{C}ech complex also bears a close relation with the derived pushforward functor ${\bm R}\pi^+_*$ between the varieties in \eqref{log_flip_intro}. Namely, every graded module $M \in \GrA$ is in particular a graded module over $A_{\geq 0}$, which therefore determines an associated sheaf%
\footnote{Notice that for each homogeneous element $f \in A$ of $\deg(f)>0$, the canonical map $(A_{\geq 0})_f \ra A_f$ of graded algebras is an isomorphism. Thus, the projective space $X^+$ is in fact covered by the spaces $\Spec(A_{(f)})$, so that we may avoid the unnatural procedure of passing to $A_{\geq 0}$.} $\widetilde{M}$ on $X^+ = \Proj(A_{\geq 0})$. Since the cohomology of quasi-coherent sheaves can be computed by \v{C}ech cohomology, we see that the weight $i$ component of $\Ce_{I^+}(M)$ may be described as
\begin{equation}  \label{CI_Rpi_intro}
{\bm R}\pi^+_* (\widetilde{M(i)}) \, \cong \, \Ce_{I^+}(M) _i
\end{equation}

This interpretation of the \v{C}ech complex will allow us to give some control over the wall-crossing problem in the case of flips and (strongly crepant) flops, in the sense of the following
\bdf  \label{flip_flop_def}
A \emph{flip} is a log flip \eqref{Xm_Y_Xp_2} in which the pair $(D^-,D^+)$ may be taken to be the canonical divisors $(K_{X^-},K_{X^+})$.

A \emph{flop}
is a log flip \eqref{Xm_Y_Xp_2} in which $K_{X^-}$ is $\bQ$-Cartier and numerically $\pi^-$-trivial; and $K_{X^+}$ is $\bQ$-Cartier and numerically $\pi^+$-trivial.

A log flip is said to be \emph{strongly crepant} if the contracted variety $Y$ (and hence $X^-$ and $X^+$) is quasi-Gorenstein, \ie, if $K_Y$ is Cartier. This implies that it is a flop.
\edf

To state our main result, we introduce some notation.
For a dualizing complex $\omega_Y^{\prime \prime \bullet}$ on $Y = \Spec \, R$, denote by $\bD_Y : \cD(\GrA)^{\op} \ra \cD(\GrA)$ the weight-degreewise dualizing functor (see \eqref{D_Y_def} below)
\begin{equation}  \label{bD_Y_intro}
\bD_Y(M)_i \, \simeq \, \RHom_{R}(M_{-i} , \omega_Y^{\prime \prime \bullet})
\end{equation}

For non-affine base $Y$, the local cohomology complex $\RGam_{I^{\pm}}(M)$ and the \v{C}ech complex $\Ce_{I^{\pm}}(M)$ sheafifies to $\RGam_{\scI^{\pm}}(\cM)$ and $\Ce_{\scI^{\pm}}(\cM)$ respectively. The dualizing functor $\bD_Y$ also extends to the non-affine case.
Denoting by $\cM(m)$ the weight shift of $\cM \in \cD(\Gr(\cA))$ by $m$, and by $\cM[n]$ the homological shift by $n$,
then our main result says
\bthm  \label{main_thm_intro}
Assume that the base field $k$ has characteristic zero. Then
\begin{enumerate}
	\item Let \eqref{log_flip_intro} be a flip, where $X^-$ and $X^+$ are projective with rational Gorenstein singularities. Take $(D^-,D^+)$ to be the canonical divisors $(K_{X^-}, K_{X^+})$, and take the corresponding sheaf $\cA := \bigoplus_{i \in \bZ} \cO_Y(i K_Y)$ of $\bZ$-graded algebras on $Y$. Then there is a dualizing complex $\omega_Y^{\prime \prime \bullet}$ on $Y$ such that the local cohomology complexes satisfy a duality $\RGam_{\scI^+}(\cA)(1)[1] \cong \bD_Y(\RGam_{\scI^-}(\cA))$. Moreover, for each open affine subscheme $U \subset Y$, the algebra $\cA(U)$ is Gorenstein.
\item Let \eqref{log_flip_intro} be a strongly crepant flop, where $X^-$ and $X^+$ are projective with rational Gorenstein singularities. Take $(D^-,D^+)$ to be Cartier, and take the corresponding sheaf $\cA := \bigoplus_{i \in \bZ} \cO_Y(i D_Y)$ of $\bZ$-graded algebras on $Y$. Then there is a dualizing complex $\omega_Y^{\prime \prime \bullet}$ on $Y$ such that the local cohomology complexes satisfy a duality $\RGam_{\scI^+}(\cA)[1] \cong \bD_Y(\RGam_{\scI^-}(\cA))$. Moreover, for each open affine subscheme $U \subset Y$, the algebra $\cA(U)$ is Gorenstein.
\end{enumerate}
\ethm

To the best of our knowledge, this is the first result that gives some control over the wall-crossing problem in general classes of flips and flops. Moreover, the information we obtain is closely related to the derived categories on two sides of the flip/flop. As we will demonstrate in Section \ref{Dcat_flip_flop_sec}, as well as in the paper \cite{Yeu20b}, this presents rather convincing evidence for the general principle \eqref{CMMP_principle_2}, and indeed a possible path to proving it.

This result is proved by formulating the Grothendieck dualities for $\pi^-$ and $\pi^+$ in terms of the \v{C}ech complex, in view of \eqref{CI_Rpi_intro}. A crucial point is a compatibility condition between the Grothendieck dualities across the flip/flop. In fact, we will see that both cases exhibit the same duality structure (which we call a homological flip/flop) depending on a paramter $a \in \bZ_{\geq 0}$, which is $a = 0$ for flops and $a = 1$ for flips.
We also remark that the conditions that $k$ has characteristic zero and that $X^{\pm}$ are projective with rational singularities are imposed in order to apply the Kawamata-Viehweg vanishing theorem%
\footnote{By using a stronger form of the Kawamata-Viehweg vanishing theorem (see Remark \ref{KV_klt_remark} below), one may weaken the assumptions of Theorem \ref{main_thm_intro}(1): it suffices to assume that $X^{\pm}$ are proper over $k$ and have log terminal singularities (\ie, $(X^{\pm},0)$ are klt). Also, in part (2), we may remove the assumption that $D^{\pm}$ are Cartier.}.
Before providing more details of this proof, we mention some applications.

First, notice that in both cases (1) and (2) of Theorem \ref{main_thm_intro}, the divisors $(D^-,D^+)$ used to define the log flip are Cartier. 
This can be used to show that Serre's equivalence holds. Indeed, in \cite{Yeu20c}, a sharpening of Serre's equivalence is established, where the usual condition that $\cA_{\geq 0}$ is generated by $\cA_1$ over $\cA_0$ is replaced by the weaker assumption that $\cA$ is ``Cartier'' (see Definition \ref{frac_Cartier} below). In the present context, this ``Cartier'' condition translates precisely to the condition that $D^{\pm}$ are Cartier (see Proposition \ref{log_flip_algebra} below), so that we have
\begin{equation*}
\QCoh(X^-) \, \simeq \, \Gr(\cA) / \Tor^-
\qquad \text{and} \qquad 
\QCoh(X^+) \, \simeq \, \Gr(\cA) / \Tor^+
\end{equation*} 
Alternatively, it says that the stacks $\mathfrak{X}^{\pm}$ are isomorphic to the varieties $X^{\pm}$ (see Lemma \ref{frac_Cartier_twist}(1)). 
As a formal consequence (see Proposition \ref{D_IpTR_DQCoh}), we have exact equivalences 
\begin{equation}  \label{D_QCoh_ITR_intro}
\cD(\QCoh(X^-)) \, \simeq \, \cD_{\scI^-\text{-triv}}(\GrA)
\qquad \text{and} \qquad
\cD(\QCoh(X^+)) \, \simeq \, \cD_{\scI^+\text{-triv}}(\GrA)
\end{equation}

By \eqref{torsion_SOD_intro}, both of \eqref{D_QCoh_ITR_intro} sit as semi-orthogonal summand of $\cD(\Gr(\cA))$. Moreover, the other components $\cD_{\Tor^{\pm}}(\Gr(\cA))$ are closely related to $\RGam_{\scI^{\pm}}(\cA)$. More precisely, the right adjoint to the inclusion $\cD_{\Tor^{\pm}}(\Gr(\cA)) \rinto \cD(\Gr(\cA))$ is given by $\cM \mapsto \cM \otimes_{\cA}^{{\bm L}} \RGam_{\scI^{\pm}}(\cA)$. 
Thus, the duality between $\RGam_{\scI^+}(\cA)$ and $\RGam_{\scI^-}(\cA)$ in Theorem \ref{main_thm_intro} should allow one to relate $\cD(\QCoh(X^-))$ and $\cD(\QCoh(X^+))$. 

An instance of this relation is given by the following Corollary. Here, we say that a full triangulated subcategory $\cE \subset \cD(\Gr(\cA))$ is $\cO_Y$-linear if $\cK \otimes_{\cO_Y}^{{\bm L}} \cM \in \cE$ for every $\cM \in \cE$ and $\cK \in \Dperf(Y)$.
\bcor  \label{Dperf_Ce_ff_flop_intro}
Let $(Y,\cA)$ be as in the situation of Theorem \ref{main_thm_intro}(2). 
Then for any $\cO_Y$-linear full triangulated subcategory $\cE \subset \Dperf(\Gr(\cA))$, the functor $\Ce_{\scI^+} : \cE \ra \cD_{\scI^+\text{-triv}}(\GrA)$ is fully faithful if and only if the functor $\Ce_{\scI^-} : \cE \ra \cD_{\scI^-\text{-triv}}(\GrA)$ is fully faithful.
\ecor

\bpf
Since $\cE$ is $\cO_Y$-linear, this statement holds if and only if it holds on every open affine subscheme $U \subset Y$, which becomes Corollary \ref{Dperf_Ce_ff_flop} below.
\epf

In recent years, many results about the change in the derived category under variations of GIT quotients are studied by singling out full triangulated subcategories $\cE^{\pm} \subset \Dbcoh(\mathfrak{X})$ whose restriction functor to  $\Dbcoh(\mathfrak{X}^{ss}(L^{\pm}))$ give exact equivalences (see, \eg, \cite{HL15, BFK19}). Moreover, one often has a direct comparison between $\cE^-$ and $\cE^+$, which therefore allows one to compare the derived categories under the wall-crossing. Viewed in this light, Corollary \ref{Dperf_Ce_ff_flop_intro} roughly says that, if we work with $\Dperf(-)$ instead of $\Dbcoh(-)$, then a full subcategory $\cE \subset \Dperf(\GrA)$ works for $X^-$ if and only if it works for $X^+$. It could therefore be viewed as an approximation for the statement \eqref{CMMP_principle_2} for flops.
In particular, we suggest in Remark \ref{Bondal_Orlov_remark} below how this might open up a new way to tackle a conjecture of Bondal and Orlov.

In the paper \cite{Yeu20b}, we develop the technique of weight truncation parallel to \cite{HL15, BFK19}, which gives another way to exploit the conclusion of Theorem \ref{main_thm_intro} to give statements along the lines of \eqref{CMMP_principle_2}.
Specifically we fall short of establishing a general statement of \eqref{CMMP_principle_2} only because certain spectral seqeuences are not known to converge (see \cite[Remark 8.7]{Yeu20b}). This convergence issue seems to be similar to phenomena one sees in Koszul duality. It seems possible that a formal modification of our argument could lead to more satisfactory statements. 

In any case, these formal problems do not arise if $\cA$ is smooth. Thus, we have the following result (see \cite[Corollary 8.6]{Yeu20b})
\bthm
\begin{enumerate}
	\item In the situation of Theorem \ref{main_thm_intro}(1), if $\cA$ is smooth over $k$, then $X^{\pm}$ are smooth%
	\footnote{Here, we use the fact that $D^{\pm}$ are Cartier to conclude that $\mathfrak{X}^{\pm} = X^{\pm}$. The Cartier condition can be weakened, in which case there is a comparison between $\Dperf(\mathfrak{X}^-)$ and $\Dperf(\mathfrak{X}^+)$.}, and there is a fully faithful exact functor $\Dbcoh(X^+) \rinto \Dbcoh(X^-)$.
	\item In the situation of Theorem \ref{main_thm_intro}(2), if $\cA$ is smooth over $k$, then $X^{\pm}$ are smooth, and there is an exact equivalence $\Dbcoh(X^+) \simeq \Dbcoh(X^-)$.
\end{enumerate}
\ethm

Similar results have been obtained in \cite{HL15, BFK19}. However, the information we use to control the wall-crossing is quite different from that in {\it loc.cit.} For example, in \cite{HL15}, to establish the analogue of (2), one requires the canonical bundle to have weight $0$ along the fixed loci. In contrast, we impose conditions on the varieties \eqref{log_flip_intro}, which mostly concerns the complement of the fixed loci. 
In practice, it is often easier to obtain information about the varieties \eqref{log_flip_intro} than about the stack $\mathfrak{X}$. 
In particular, the conditions in Theorem \ref{main_thm_intro}(2) are almost automatically satisfied if $X^{\pm}$ are projective Calabi-Yau varieties with rational Gorenstein singularities. 

In the rest of this introduction, we explain our proof of Theorem \ref{main_thm_intro}. First, we specify the choice of the dualizing complexes in the Theorem. Let $\pi : Y \ra \Spec \, k$ be the projection to a point, and let $\omega_Y^{\bullet} := \pi^!(k)$ be the canonical dualizing complex on $Y$. Let $n = \dim Y$.
\begin{enumerate}
	\item In the setting of Theorem \ref{main_thm_intro}(1), let $\omega_Y^{\prime \prime \bullet} := \omega_Y^{\bullet}[-n]$.
	\item In the setting of Theorem \ref{main_thm_intro}(2), let $\omega_Y^{\prime \prime \bullet} := \omega_Y^{\bullet}[-n] \otimes_{\cO_Y} \cO(-K_Y)$.
\end{enumerate}

We write $\cO_{X^{\pm}}(i) := \cO_{X^{\pm}}(i D^{\pm})$. Then in both case, there are isomorphisms
\begin{equation}  \label{pi_shriek_omega_intro}
(\pi^{\pm})^{!}(\omega_Y^{\prime \prime \bullet}) \, \cong \, \cO_{X^{\pm}}(a)
\end{equation}
where $a = 1$ in the setting of Theorem \ref{main_thm_intro}(1), and $a = 0$ in the setting of Theorem \ref{main_thm_intro}(2).

Under the isomorphisms \eqref{pi_shriek_omega_intro}, the local adjunction isomorphisms then read
\begin{equation}  \label{local_adj_flop_intro}
{\bm R}\pi^{\pm}_* \cO_{X^{\pm}}(a+i) \, = \, {\bm R}\pi^{\pm}_* \RcHom_{\cO_{X^{\pm}}}( \cO_{X^{\pm}}(-i) , \cO_{X^{\pm}}(a) ) \, \cong \, \RcHom_{\cO_{Y}}( {\bm R}\pi^{\pm}_*\cO_{X^{\pm}}(-i) , \omega_Y^{\prime \prime \bullet} )
\end{equation}

In view of \eqref{CI_Rpi_intro}, the left-hand-side is the weight $i$ component of $\Ce_{\scI^{\pm}}(\cA)(a)$, while the right-hand-side is the weight $i$ component of $\bD_Y( \Ce_{\scI^{\pm}}(\cA) )$. This suggests that there are isomorphisms in $\cD(\Gr(\cA))$:
\begin{equation}  \label{Phi_pm_intro}
\Phi^{\pm} \, : \, \Ce_{\scI^{\pm}}(\cA)(a) \raq \bD_Y(\Ce_{\scI^{\pm}}(\cA))
\end{equation}


The isomorphisms \eqref{Phi_pm_intro}, per se, are not so useful as each of them concerns only one side of the flip/flop. In order to relate the two sides of the flip/flop, a crucial observation is that $\Phi^-$ and $\Phi^+$ satisfy a certain compatibility condition, as expressed by the commutativity of the diagram \eqref{HFF_compatible_intro} below.

\bdf  \label{HFF_def_intro}
A \emph{homological flip} (resp. \emph{homological flop}) consists of a sextuple $(Y,\omega_Y^{\prime \prime \bullet}, \cA, a, \Phi^-,\Phi^+)$ where
\begin{enumerate}
	\item $Y$ is a Noetherian separated scheme with a dualizing complex $\omega_Y^{\prime \prime \bullet} \in \Dbcoh(\QCoh(Y))$.
	\item $\cA$ is a quasi-coherent sheaf of Noetherian $\bZ$-graded rings such that $\cA_0$ is coherent over $Y$.
	\item $a > 0$ (resp. $a=0$) is an integer
	\item $\Phi^-$ and $\Phi^+$ are isomorphisms in $\cD(\Gr(\cA))$ 
	\begin{equation*}  
	\begin{split}
	\Phi^+ \, &: \, \Ce_{\scI^+}(\cA)(a) \raq \bD_Y(\Ce_{\scI^+}(\cA)) \\
	\Phi^- \, &: \, \Ce_{\scI^-}(\cA)(a) \raq \bD_Y(\Ce_{\scI^-}(\cA))
	\end{split}
	\end{equation*}
	where $\bD_Y : \cD(\Gr(\cA))^{\op} \ra \cD(\Gr(\cA))$ is defined using the dualizing complex $\omega_Y^{\prime \prime \bullet}$, as in \eqref{bD_Y_intro}.
\end{enumerate}

The maps $\Phi^-$ and $\Phi^+$ are required to be compatible in the sense that the diagram 
\begin{equation}  \label{HFF_compatible_intro}
\begin{tikzcd} [row sep = 0]
& \Ce_{\scI^+}(\cA)(a) \ar[r, "\Phi^+"] & \bD_Y( \Ce_{\scI^+}(\cA) ) \ar[rd, "\bD_Y(\eta^+)"] & \\
\cA(a) \ar[ru, "\eta^+"] \ar[rd, "\eta^-"']  & & & \bD_Y(\cA) \\
& \Ce_{\scI^-}(\cA)(a) \ar[r, "\Phi^-"] & \bD_Y( \Ce_{\scI^-}(\cA) ) \ar[ru, "\bD_Y(\eta^-)"'] &
\end{tikzcd}
\end{equation}
commutes in $\cD(\Gr(\cA))$, where $\eta^{\pm} : \cA \ra \Ce_{\scI^{\pm}}(\cA)$ are the maps in \eqref{RGI_CI_seq_intro}.
\edf

Our above discussion suggests that the settings of Theorem \ref{main_thm_intro} constitutes an example of homological flip/flop. This is indeed true, and is proved in Theorem \ref{HFF_main_examples} below, which is formulated slightly more generally than Theorem \ref{main_thm_intro}.

In order to establish the compatibility condition \eqref{HFF_compatible_intro}, one needs to unravel the definitions of the isomorphisms  \eqref{CI_Rpi_intro} and \eqref{pi_shriek_omega_intro}, and see how they translate to \eqref{local_adj_flop_intro}, and then to \eqref{Phi_pm_intro}. For these purposes, we have found it necessary to establish from scratch some of the more or less known results concerning Serre's equivalence and Grothendieck duality. This is carried out in details in \cite{Yeu20c}.
The results we need are summarized in Section \ref{DGrA_sec} below.
On the other hand, the isomorphism \eqref{pi_shriek_omega_intro} follows from a standard result that $\cO(K_X) \cong \omega_X$ for normal projective varieties $X$. A proof of this isomorphism is given in \cite{KM98}, with a second proof being hinted at.
We provide some details of this second proof in Section \ref{flip_flop_subsec} below (see Theorem \ref{O_KX_dualizing}) in order to keep track of some compatibility conditions (Lemma \ref{compatible_lemma_1}, \ref{compatible_lemma_2}).

One can say more in the situation of Theorem \ref{main_thm_intro}. Namely, by the Kawamata-Viehweg vanishing theorem (see Corollary \ref{KV_vanishing_2}), one has ${\bm R}^j \pi^{+}_*  (\cO_{X^+}(i)) = 0$ for $i \geq a$ and ${\bm R}^j \pi^{-}_* ( \cO_{X^-}(i)) = 0$ for $i \leq a$.
This translates to the following condition:

\bdf  \label{canonical_vanishing_intro}
A homological flip/flop $(Y,\omega_Y^{\prime \prime \bullet}, \cA, a, \Phi^-,\Phi^+)$ is said to have \emph{canonical vanishing} if $\RGam_{\scI^+}(\cA)_i \simeq 0$ for all $i \geq a$ and $\RGam_{\scI^-}(\cA)_i \simeq 0$ for all $i \leq a$.
\edf

Our Theorem \ref{main_thm_intro} then follows from the following
\bthm  \label{HFF_Gorenstein}
Let $(Y,\omega_Y^{\prime \prime \bullet}, \cA, a, \Phi^-,\Phi^+)$ be a homological flip/flop with canonical vanishing, then there is a map $\Psi : \RGam_{\scI^+}(\cA)(a)[1] \xra{\cong} \bD_Y(\RGam_{\scI^-}(\cA))$ in $\cD(\Gr(\cA))$ that is a quasi-isomorphism in all weights $i \in \bZ$ except possibly $-a < i < 0$.

In particular, if $a = 0$ or $a=1$, then $\Psi$ is an isomorphism. As a consequence, for each open affine $U \subset Y$, the ring $\cA(U)$ is Gorenstein.
\ethm

This result follows from Theorem \ref{HFF_main_thm} below (see the end of Section \ref{HFF_sec}). Notice that after we have extracted the homological structure of the flip/flop in the form of Definition \ref{HFF_def_intro} and \ref{canonical_vanishing_intro}, the proofs of Theorem \ref{HFF_main_thm} and \ref{HFF_Gorenstein} proceed  by pure homological/commutative algebra, and is therefore self-contained. 

%

\vspace{0.3cm}

\textbf{Acknowledgements.} The author would like to thank Yuri Berest, Daniel Halpern-Leistner, and Valery Lunts for helpful discussions.

\vspace{0.3cm}

\textbf{Conventions.} A variety means an integral separated scheme of finite type over a field. We will always work with cochain complexes. We denote the category of cochain complexes of an abelian category $\cC$ by $\Ch(\cC)$. An inverse functor will mean what is sometimes called a quasi-inverse (\ie, composition on both sides are isomorphic to the identity functors).

\section{Derived categories of graded modules}  \label{DGrA_sec}

In this section, we review some basic results about graded modules over graded rings.
See \cite{Yeu20c} for unexplained notation and terminology. 

\bdf
A \emph{$\bZ$-graded ring} is a commutative ring $A$ with a $\bZ$-grading $A = \bigoplus_{n \in \bZ} A_n$. 
Here, by commutative we mean $xy = yx$, not $xy = (-1)^{|x||y|} yx$. 

A \emph{graded module} over $A$ will always mean a $\bZ$-graded module $M = \bigoplus_{n \in \bZ} M_n$.
\edf

We first recall the following result (see, e.g., \cite[Theorem 1.5.5]{BH93}):
\bpp  \label{Noeth_gr_ring}
Let $A$ be a $\bZ$-graded ring. Then the followings are equivalent:
\begin{enumerate}
	\item $A$ is a Noetherian ring; 
	\item every graded ideal of $A$ is finitely generated;
	\item $A_0$ is Noetherian, and both $A_{\geq 0}$ and $A_{\leq 0}$ are finitely generated over $A_0$;
	\item $A_0$ is Noetherian, and $A$ is finitely generated over $A_0$.
\end{enumerate}
\epp

Denote by $\GrA$ the category of graded modules over $A$, whose morphisms are maps of graded modules of degree $0$.
Denote by $\Homcom_A(M,N)$ the graded module of Hom complexes, so that $\Hom_{\GrA}(M,N) =  \Hom_A(M,N) = \Homcom_A(M,N)_0$.
Let $\RHomcom_A(M,N)$ be the derived version, then we have the following graded analogue of 
\cite[Tag 0ATK]{Sta}:

\bpp  \label{tensor_in_Hom_target}
Assume that $A$ is Noetherian, then for any $N \in \Dmcoh(\GrA)$, $L \in \cD^+(\GrA)$, and $M \in \cD(\GrA)$ of finite Tor dimension, the canonical map
\begin{equation*}
M \otimes_A^{{\bm L}} \RHomcom_A(N,L) \raq \RHomcom_A(N, M\otimes_A^{{\bm L}} L)
\end{equation*}
is an isomorphism in $\cD(\GrA)$.
\epp

Let $I \subset A$ be a finitely generated graded ideal, then there is a semi-orthogonal decomposition (see \cite[(2.41)]{Yeu20c}) of the derived category $\cD(\GrA)$ of graded $A$-modules
\begin{equation}  \label{local_cohom_SOD}
\cD(\GrA) \, = \, \langle \, \cD_{\ITR}(\GrA) \, , \,  \cD_{\IIT}(\GrA) \, \rangle 
\end{equation}
such that for every $M \in \cD(
\GrA)$, the corresponding decomposition sequence is the local cohomology sequence 
\begin{equation}  \label{RGam_Ce_seq}
\ldots \raq \RGI(M) \xraq{\epsilon_M} M \xraq{\eta_M} \CI(M) \xraq{\delta_M} \RGI(M)[1] \raq \ldots 
\end{equation}
which is an obvious generalization of \eqref{RGI_CI_seq_intro} from the case $M \in \GrA$ to $M \in \cD(\GrA)$.

In fact, the semi-orthogonal decomposition \eqref{local_cohom_SOD} is only part of a recollement (see \cite{Yeu20c})
\begin{equation}  \label{Greenlees_May_rec}
\begin{tikzcd} [row sep = 0, column sep = 50]
& & \cD_{\IIT}(\GrA) \ar[ld, bend right, "\iota"'] \ar[dd, shift left = 3, "\LLI"]  \ar[dd, phantom, "\simeq"] \\
\cD_{\ITR}(\GrA) \ar[r, "\iota"] 
& \cD(\GrA) \ar[ru, "\RGI"] \ar[rd, "\LLI"']  \ar[l, bend right, "\Ce_I"'] \ar[l, bend left, "\cE_I"]  \\
& & \cD_{\Icomp}(\GrA)  \ar[lu, bend left, "\iota"] \ar[uu, shift left = 3 , "\RGI"]
\end{tikzcd}
\end{equation}

We are mostly interested in the special case $I = I^+ := A_{>0} \cdot A$, where the corresponding recollement \eqref{Greenlees_May_rec} is closely related to the projective space $X^+ = \Proj^+(A) := \Proj(A_{\geq 0})$.
Indeed, if we denote by $\pi^+ : \Proj^+(A) \ra \Spec(A_0) =: Y$ the canonical map, then there are three aspects of this relation:
\begin{equation}  \label{GGM_three_points}
\parbox{40em}{(1) The category $\cD_{\IpTR}(\GrA)$ is closely related to the category $\cD(\QCoh(X^+))$.\\
(2) The functor $\Ce_{I^+}$ is closely related to the derived pushforward functor ${\bm R}\pi^+_*$.\\
(3) The functor $\cE_{I^+}$ is closely related to the shriek pullback functor $(\pi^+)^!$.}
\end{equation}

Point (2) is easiest to see. Namely, for each $M \in \cD(\GrA)$ and for each $i \in \bZ$, there is a canonical isomorphism in $\cD(A_0)$:
\begin{equation}  \label{Ce_i_RGam}
\Ce_{I^+}(M)_i \, \cong \, {\bm R}\pi^+_*\,  \widetilde{M(i)} 
\end{equation}
where we denote by $\widetilde{M} \in \QCoh(X^+)$ the quasi-coherent sheaf on $X^+$ associated to $M \in \GrA$.

This gives a proof of the following two results (see \cite[Lemma 4.12, 4.13]{Yeu20c}):
\blm  \label{RGam_Dbcoh}
If $A$ is Noetherian, then for any $M \in \Dbcoh(\GrA)$, we have $\Ce_{I^+}(M)_i \in \Dbcoh(A_0)$ and $\RGam_{I^+}(M)_i \in \Dbcoh(A_0)$ for each weight $i \in \bZ$.
\elm

\blm  \label{local_cohom_weight_bounded}
If $A$ is Noetherian, then for any $M \in \Dbcoh(\GrA)$, there exists $c^+, c^- \in \bZ$ such that 
$\RGam_{I^+}(M)_i \simeq 0$ for all $i \geq c^+$ and $\RGam_{I^-}(M)_i \simeq 0$ for all $i \leq c^-$.
\elm

Point (1) of \eqref{GGM_three_points} concerns Serre's equivalence. Our version stresses the following condition:

\bdf  \label{frac_Cartier}
A $\bZ$-graded ring $A$ is said to be \emph{positively $\tfrac{1}{d}$-Cartier}, for an integer $d > 0$, if 
the canonical map $\widetilde{A(di)} \otimes_{\cO_{X^+}} \widetilde{A(dj)} \ra \widetilde{A(di+dj)}$ is an isomorphism for all $i,j \in \bZ$.
Similarly, it is said to be \emph{negatively $\tfrac{1}{d}$-Cartier} if the analogous condition holds for $X^- = \Proj^-(A)$.
It is said to be \emph{$\tfrac{1}{d}$-Cartier} if it is both positively and negatively $\tfrac{1}{d}$-Cartier.
When $d=1$, we simply say that $A$ is \emph{(positively/negatievly) Cartier}.
\edf

For example, if  $A_{\geq 0}$ is generated over $A_0$ by homogeneous elements $f_1,\ldots,f_p$ of positive degrees $d_i := \deg(f_i) > 0$, then it can be shown that $A$ is positively $\tfrac{1}{d}$-Cartier for any $d > 0$ that is divisible by each of $d_i$ (see \cite[Lemma 3.5]{Yeu20c}).
However, this sufficient condition is by no means necessary. In examples arising from log flips, the Cartier condition translates to certain divisors being Cartier (see Proposition \ref{log_flip_algebra} below). This is more natural and easier to check than the condition that $A_{\geq 0}$ be generated by $A_1$ over $A_0$.

\blm  \label{frac_Cartier_twist}
If $A$ is positively $\tfrac{1}{d}$-Cartier, then 
\begin{enumerate}
	\item $\Proj^+(A)$ can be covered by the standard open subsets $D(fg)$, where $f , g$ are homogeneous elements of positive degrees such that  $\deg(f) - \deg(g) = d$.
	\item $\widetilde{A(di)}$ is an invertible sheaf for all $i \in \bZ$, and the canonical map $\widetilde{M} \otimes_{\cO_{X^+}} \widetilde{A(di)} \ra \widetilde{M(di)}$ is an isomorphism for any $M \in \GrA$.
\end{enumerate}
\elm

One version of Serre's equivalence takes following form (see \cite[Theorem 3.15]{Yeu20c}):

\bthm  \label{Serre_equiv_v1}
Suppose that $A$ is Noetherian and positively Cartier, then there is an equivalence of categories 
$\QCoh(X^+) \simeq \QpGr(A) := \GrA / \Torp$.
\ethm

A more detailed discussion of Serre's equivalence will be delegated to the end of this section, where we work with the setting \eqref{sheaf_A_setting} of non-affine base, and is therefore notationally more convenient for the rest of this paper.

Point (3) of \eqref{GGM_three_points} concerns an interpretation of Grothendieck duality in terms of Greenlees-May duality%
\footnote{By Greenlees-May duality, we mean the recollement \eqref{Greenlees_May_rec}.} of graded modules. 
Assume that $A$ is Noetherian, then the adjunction isomorphism between $\widetilde{A(i)} \in \Coh(X^+)$ and any $L \in \QCoh(Y) \simeq \Mod(A_0)$ gives an isomorphism in $\cD(A_0)$:
\begin{equation}  \label{adj_Ai_L_1} 
\RHom_{\cO_{X^+}}( \widetilde{A(-i)} , (\pi^+)^! L)
\, \cong \, 
\RHom_{A_0}( {\bm R} \pi^+_*( \widetilde{A(-i)} ) , L)
\end{equation}

Both sides of this isomorphism are naturally the weight $i$ component of objects in $\cD(\GrA)$. For the left hand side, consider the right adjoint $\cR^+ : \cD(\QCoh(X^+)) \ra \cD(\GrA)$ to the associated sheaf functor $(-)^{\sim}$. One can show that the left hand side of \eqref{adj_Ai_L_1} is the weight $i$ component of $\cR^+((\pi^+)^!(L))$. On the other hand, by \eqref{Ce_i_RGam}, the right hand side of \eqref{adj_Ai_L_1} is the weight $i$ component of $\RHomcom_{A_0}(\Ce_{I^+}(A), L)$.
This suggests the following
\bthm[\cite{Yeu20c}, Theorem 4.16]  \label{GGM_thm}
For each $L \in \cD(A_0)$, there is a canonical isomorphism in $\cD(\GrA)$:
\begin{equation*}
\cR^+ ((\pi^+)^! ( L ))  \, \cong \,  \cE_{I^+}( \RHomcom_{A_0}(A,L) ) 
\, \cong \, \RHomcom_{A_0}(\Ce_{I^+}(A),L)
\end{equation*}
whose weight $i$ component is the adjunction isomorphism \eqref{adj_Ai_L_1} under the above canonical isomorphisms. 
In other words, the isomorphisms \eqref{adj_Ai_L_1}, for $i \in \bZ$, are organized into an isomorphism in $\cD(\GrA)$.
\ethm

%

Most of these results, as summarized in \eqref{Greenlees_May_rec} and \eqref{GGM_three_points}, carries over to the case of non-affine base. More precisely, we consider the following setting:
\begin{equation}  \label{sheaf_A_setting}
\parbox{40em}{$Y$ is a Noetherian separated scheme, and $\cA$ is a quasi-coherent sheaf of Noetherian $\bZ$-graded rings on $Y$, such that $\cA_0$ (and hence every $\cA_i$) is coherent over $\cO_Y$.}
\end{equation}

In this setting, we consider the category $\Gr(\cA)$ of quasi-coherent graded $\cA$-modules, as well as its derived category $\cD(\Gr(\cA))$. For each quasi-coherent sheaf of graded ideals $\scI \subset \cA$, the semi-orthogonal decomposition \eqref{local_cohom_SOD}, as well as the local cohomology sequence \eqref{RGam_Ce_seq}, carries over without change (see \cite[Section 5]{Yeu20c} for details).
The same is true for (1) and (2) of \eqref{GGM_three_points}. 

For example, for point (2) of \eqref{GGM_three_points},
consider $X^+ := \Proj_Y^+(\cA)$, with the canonical map $\pi^+ : X^+ \ra Y$. Then for any $\cM \in \cD(\Gr(\cA))$, we again have a canonical isomorphism 
\begin{equation}  \label{Ce_i_RGam_sheaf}
\Ce_{\scI^+}(\cM)_i \, \cong \, {\bm R}\pi^+_*\,  \widetilde{\cM(i)} 
\end{equation}
In particular, for $\cM \in \Gr(\cA)$, this gives a description of the weight components of the map $\eta_{\cM} : \cM \ra \Ce_{\scI^+}(\cM)$ in \eqref{RGam_Ce_seq}. Namely, the following diagram in $\cD(\QCoh(Y))$ commutes (see \cite[(5.37)]{Yeu20c}):
\begin{equation}  \label{eta_M_weight_description}
\begin{tikzcd}
\cM_i \ar[rr, "\eta_{\cM}"] \ar[d, equal] & &  \Ce_{\scI^+}(\cM)_i \ar[d, "\eqref{Ce_i_RGam_sheaf}", "\cong"'] \\
\cM_i \ar[r, "\eqref{M_0_pi_M}"]  &  \pi^+_*(\widetilde{\cM(i)}) \ar[r, "\iota"]  &   {\bm R}\pi^+_*(\widetilde{\cM(i)})
\end{tikzcd}
\end{equation}

However, there is some technical subtleties when one tries to define the sheafified version of local homology (\ie, the analogue of the functor $\LLI$ in \eqref{Greenlees_May_rec}). 
Indeed, already in the ungraded case, the notion of local homology has a technical subtlety arising from the fact that the internal Hom between quasi-coherent sheaves may not be quasi-coherent (see, e.g., \cite[Remark (0.4)]{ATJLJ97} for a discussion). The notion of (quasi-)coherator from \cite[Appendix B]{TT90} is therefore highly relevant in this discussion. We recall this notion now.

On a quasi-compact separated scheme $X$, the inclusion $\iota : \QCoh(X) \ra \Mod(\cO_X)$ has a right adjoint, called the \emph{(quasi-)coherator} $Q_X : \Mod(\cO_X) \ra \QCoh(X)$. The derived functor ${\bm R}Q_X : \cD(X) \ra \cD(\QCoh(X))$ then restricts to an equivalence ${\bm R}Q_X : \Dqcoh(X) \ra \cD(\QCoh(X))$, which is inverse to the inclusion
(see, e.g., \cite[Tag 08DB]{Sta}). 
Throughout this paper, we work with $\cD(\QCoh(X))$ instead of $\Dqcoh(X)$. More precisely, we work with the following two conventions:
\begin{conv}  \label{Dcat_conv1}
	On a quasi-compact separated scheme $X$, let $\cD(X) = \cD(\Mod(\cO_X))$ be the derived category of all sheaves of $\cO_X$-modules.
	
	For any $\cF,\cG \in \Mod(\cO_X)$, denote by $\cHom^{\clubsuit}_{\cO_X}(\cF,\cG) \in \Mod(\cO_X)$ the 
	internal Hom object in $\Mod(\cO_X)$, defined by $\cHom^{\clubsuit}_{\cO_X}(\cF,\cG)(U) = \Hom_{\cO_U}(\cF|_U, \cG|_U)$.
	
	For any $\cF,\cG \in \QCoh(X)$, denote by $\cHom_{\cO_X}(\cF,\cG) \in \QCoh(X)$ the 
	internal Hom object in $\QCoh(X)$, defined by $\cHom_{\cO_X}(\cF,\cG) = Q_X( \cHom^{\clubsuit}_{\cO_X}(\cF,\cG) )$.
	
	Both of these internal Hom functors have derived functors
	\begin{equation*}
		\begin{split}
			\RcHom^{\clubsuit}_{\cO_X}(-,-) \, &: \, \cD(X)^{\op} \times \cD(X) \raq \cD(X) \\
			\RcHom_{\cO_X}(-,-) \, &: \, \cD(\QCoh(X))^{\op} \times \cD(\QCoh(X)) \raq \cD(\QCoh(X))
		\end{split}
	\end{equation*}
	so that for any $\cF,\cG \in \cD(\QCoh(X))$, we have 
	\begin{equation*}
		\RcHom_{\cO_X}(\cF,\cG) \, \cong \, {\bm R}Q_X( \RcHom^{\clubsuit}_{\cO_X}(\cF,\cG) )
	\end{equation*}
	
	In particular, if $X$ is Noetherian, $\cF \in \Dmcoh(\QCoh(X))$ and $\cG \in \cD^+(\QCoh(X))$, then the canonical map 
	$ \RcHom_{\cO_X}(\cF,\cG) \ra  \RcHom^{\clubsuit}_{\cO_X}(\cF,\cG)$ is an isomorphism in $\cD(X)$ (see, e.g., \cite[Tag 0A6H]{Sta}). As a result, these two will often be implicitly identified in this case.
\end{conv}

\begin{conv}  \label{Dcat_conv2}
	Let $f : X \ra Y$ be a quasi-compact morphism between quasi-compact separated schemes. The functors $({\bm L}f^*, {\bm R}f_*, f^!)$ are regarded as functors between $\cD(\QCoh(X))$ and $\cD(\QCoh(Y))$. Namely, ${\bm L}f^*$ and ${\bm R}f_*$ are the derived functors of the functors $f^*$ and $f_*$ between $\QCoh(X)$ and $\QCoh(Y)$.
	Under the equivalence ${\bm R}Q : \Dqcoh(-) \xra{\simeq} \cD(\QCoh(-))$, the functor ${\bm R}f_*$ coincides with the usual one between $\Dqcoh(X)$ and $\Dqcoh(Y)$ (see, e.g., \cite[Tag 0CRX]{Sta}). By adjunction, the same holds for ${\bm L}f^*$, as well as $f^!$, whenever it is well-defined (see Section \ref{flip_flop_subsec} for a summary of the functor $f^!$).
\end{conv}

There is also a (quasi-)coherator for sheaves of $\cA$-modules, which is in fact computed by taking the (quasi-)coherator in each weight component (see, \eg, the proof of \cite[Proposition 5.16]{Yeu20c}). One can then use it to define the internal graded Hom complex in $\cD(\Gr(\cA))$ (see \cite[Section 5]{Yeu20c} for details):
\begin{equation*}
\RcHomcom_{\cA}(\cM,\cN) \, := \,  {\bm R}Q_{\cA} \, \RcHomcom^{\clubsuit}_{\cA}(\cM,\cN)
\end{equation*}
This Hom complex is local on $Y$ if $\cM \in \Dmcoh(\Gr(\cA))$ and $\cN \in \cD^+(\Gr(\cA))$ (see \cite[Remark 5.17]{Yeu20c}).

We will also have occasions to use a degreewise sheafified Hom complex. Namely, for any $\cM \in \cD(\Gr(\cA))$ and $\cG \in \cD(\QCoh(Y))$, there is a naturally defined Hom object $\RcHomcom_{\cO_Y}(\cM , \cG) \in \cD(\Gr(\cA))$ whose weight $i$ component is given by
\begin{equation*}
\RcHomcom_{\cO_Y}(\cM , \cG)_i \, \cong \, \RcHom_{\cO_Y}(\cM_{-i} , \cG)
\end{equation*}
This Hom complex is local on $Y$ if $\cM_i \in \Dmcoh(\QCoh(Y))$ for each $i \in \bZ$, and $\cN \in \cD^+(\QCoh(Y))$ (see \cite[Lemma 5.22]{Yeu20c}).

These two Hom complexes satisfies the adjunction
\begin{equation}  \label{RHom_A_Y_adj}
\RcHomcom_{\cO_Y}(\cM \otimes_{\cA}^{{\bm L}} \cN, \cG)
\, \cong \, 
\RcHomcom_{\cA}(\cM, \RcHomcom_{\cO_Y}(\cN,\cG))
\end{equation}

An important special case is when $\cG$ is a dualizing complex $\omega^{\prime \prime \bullet}_Y \in \Dbcoh(\QCoh(Y))$. In this case, $\RcHomcom_{\cO_Y}(\cM,\omega^{\prime \prime \bullet}_Y)$ is local on the base when $\cM \in \cD_{\lc}(\Gr(\cA))$, in the sense of the following 
\bdf
A quasi-coherent sheaf $\cM \in \Gr(\cA)$ of graded $\cA$-modules is said to be \emph{locally coherent} if each $\cM_i$ is coherent as a sheaf over $Y$. 
For each $\spadesuit \in \{ \, \, , +,-,b\}$, denote by $\Dsuit_{\lc}(\Gr(\cA))$ the full subcategory of $\Dsuit(\Gr(\cA))$ consisting of objects with locally coherent cohomology sheaves.
\edf

Given a dualizing complex $\omega^{\prime \prime \bullet}_Y \in \Dbcoh(\QCoh(Y))$, consider the functor
\begin{equation}  \label{D_Y_def}
\bD_Y \, : \, \cD_{\lc}(\Gr(\cA))^{\op} \raq \cD_{\lc}(\Gr(\cA)) \, , \qquad \cM \mapsto \RcHomcom_{\cO_Y}(\cM,\omega_Y^{\prime \prime \bullet})
\end{equation}
which is an involutive equivalence.

With these conventions about Hom complexes understood, if we denote by $\cR^+ : \cD(\QCoh(X^+)) \ra \cD(\Gr(\cA))$ the right adjoint of the associated sheaf functor, then we again have
\begin{equation}  \label{Rp_sheaf}
\cR^+(\cF)_i \, \cong \, {\bm R}\pi^+_*\, \RcHom_{\cO_{X^+}}(\widetilde{A(-i)}, \cF) 
\end{equation}

There is a canonical map (see \cite[(5.40)]{Yeu20c})
\begin{equation}  \label{Cech_R_functor_nonaffine}
\Ce_{\scI^+}(\cM) \raq \cR^+(\widetilde{\cM})
\end{equation}
whose weight $i$ component is described by the following commutative diagram in $\cD(\QCoh(Y))$:
\begin{equation} \label{Cech_R_functor_nonaffine_comp}
\begin{tikzcd}
\Ce_{\scI^+}(\cM)_i \ar[r, "\eqref{Cech_R_functor_nonaffine}"] \ar[d, equal, "\eqref{Ce_i_RGam_sheaf}"']
& \cR^+(\widetilde{\cM})_i \ar[d, equal, "\eqref{Rp_sheaf}"]\\
{\bm R}\pi^+_* ( \widetilde{\cM(i)})  \ar[r]
& {\bm R}\pi^+_* \RcHom_{\cO_{X^+}}(\widetilde{\cA(-i)}, \widetilde{\cM})
\end{tikzcd}
\end{equation}
where the bottom map is the canonical one.

There is also the following analogue of Theorem \ref{GGM_thm}, written in a more precise form:

\bthm [\cite{Yeu20c}, Theorem 5.42]  \label{GGM_thm_nonaffine}
For each $\cG \in \cD(\QCoh(Y))$, there is a canonical isomorphism in $\cD(\Gr(\cA))$: 
\begin{equation}  \label{cR_pi_shriek_isom_sheaf}
\cR^+ ((\pi^+)^! ( \cG ))  \, \cong \,  \cE_{\scI^+}( \RcHomcom_{\cO_Y}(\cA, \cG ) ) \, \cong \, 
\RcHomcom_{\cO_Y}(\Ce_{\scI^+}(\cA), \cG )
\end{equation}
Moreover, under the isomorphisms \eqref{Ce_i_RGam_sheaf} and \eqref{Rp_sheaf}, the weight $i$ component of this isomorphism can be described by the commutativity of the following diagram in $\cD(\QCoh(Y))$:
\begin{equation}  \label{cR_pi_shriek_isom_2}
\begin{tikzcd}
\cR^+ ((\pi^+)^! ( \cG ))_i \ar[r, "\eqref{cR_pi_shriek_isom_sheaf}", "\cong"'] \ar[d, "\eqref{Rp_sheaf}"', "\cong"]  &   \RcHom_{\cO_Y}(\Ce_{\scI^+}(\cA)_{-i},\cG) \ar[d, "\eqref{Ce_i_RGam_sheaf}", "\cong"'] \\
{\bm R}\pi^+_* \RcHom_{\cO_{X^+}}( \widetilde{\cA(-i)} , (\pi^+)^!(\cG) ) \ar[r, "\cong"] & \RcHom_{\cO_Y}( {\bm R}\pi^+_*(\widetilde{\cA(-i)}),\cG)
\end{tikzcd}
\end{equation}
where the horizontal map in the bottom is the local adjunction isomorphism.
\ethm

\brm  \label{weight_comp_assemble_remark}
The existence of the map \eqref{Cech_R_functor_nonaffine} and the isomorphism \eqref{cR_pi_shriek_isom_sheaf} in $\cD(\GrA)$ should be intuitively quite clear, as there are canonical maps (resp. isomorphisms) relating each of their weight components. However, to actually prove that these assemble to maps in $\cD(\GrA)$, we have found it necessary to develop the relevant functors from scratch, as we do in \cite{Yeu20c}.
\erm

Now we discuss Serre's equivalence in the setting \eqref{sheaf_A_setting}. 
We define the functor 
\begin{equation}
\,^0 \! \cL^+ : \QCoh(X^+) \ra \Gr(\cA) \, , \qquad \,^0 \! \cL^+ (\cF)_i \, := \, \pi^+_*( \cF \otimes_{\cO_{X^+}} \widetilde{\cA(i)})
\end{equation}
To prove an analogue of Theorem \ref{Serre_equiv_v1}, there are two aspects to keep track of:
\begin{equation}  \label{tilde_L_relations}
\parbox{40em}{(1) For $\cM \in \Gr(\cA)$, how are $\cM$ and $\,^0 \! \cL^+(\widetilde{\cM})$ related?\\
(2) For $\cF \in \QCoh(X^+)$, how are $\cF$ and $\widetilde{\,^0 \! \cL^+(\cF)}$ related?
}
\end{equation}

For point (1), notice that for each $\cM \in \Gr(\cA)$, there is a canonical map
\begin{equation}  \label{M_0_pi_M}
\cM_0 \raq \pi^+_*(\widetilde{\cM})
\end{equation}
Thus, for each $i \in \bZ$, there are canonical maps
\begin{equation*}  
\cM_i \raq \pi^+_*(\widetilde{\cM(i)}) \laq \pi^+_*(\widetilde{\cM} \otimes_{\cO_{X^+}} \widetilde{\cA(i)}) 
\end{equation*}
which assemble to give maps of graded modules, \ie, they are the weight $i$ components of the maps 
\begin{equation}  \label{M_Ce_L}
\cM \raq \,^0 \Ce_{\scI^+}(\cM) \laq \,^0 \! \cL^+(\widetilde{\cM})
\end{equation}

The notion of a Cartier $\bZ$-graded ring in Definition \ref{frac_Cartier} extends directly to the present setting of non-affine base. The following proposition answers the question \eqref{tilde_L_relations}(1):

\bpp  \label{M_Ce_L_isoms}
If $\cA$ is Noetherian and positively $\tfrac{1}{d}$-Cartier, then 
\begin{enumerate}
\item Both the kernel and cokernel of the map $\cM \ra \,^0 \Ce_{\scI^+}(\cM)$ in \eqref{M_Ce_L} are $(\scI^+)^{\infty}$-torsion.
\item The map $\,^0 \! \cL^+(\widetilde{\cM}) \ra \,^0 \Ce_{\scI^+}(\cM)$ in \eqref{M_Ce_L} is an isomorphism in weight $i \in d\bZ$.
\end{enumerate}
\epp

For the question \eqref{tilde_L_relations}(2), we consider the following more general setting:
\begin{equation}  \label{sheaf_of_alg_on_X}
\parbox{40em}{
	On a scheme $X$ proper over another Noetherian scheme $Y$, there is a quasi-coherent sheaf of $\bZ$-graded algebras $\bigoplus_{i \in \bZ} \cO(i)$ such that each $\cO(i)$ is coherent, and there exists some positive integer $d>0$ such that\\
	(1) the sheaf $\cO(d)$ is an invertible sheaf on $X$ ample over $Y$; \\
	(2) for each $i,j \in \bZ$, the multiplication map $\cO(di) \otimes_{\cO_{X}} \cO(j) \ra \cO(di+j)$ is an isomorphism.}
\end{equation}

For example, if we are given a positively $\tfrac{1}{d}$-Cartier pair $(Y,\cA)$ as in  \eqref{sheaf_A_setting}, then by Lemma \ref{frac_Cartier_twist}, $X := X^+$ and $\cO(i) := \widetilde{\cA(i)}$ satisfies this setting. However, for the examples of flips and flops that we consider in this paper, we actually start with \eqref{sheaf_of_alg_on_X} and consider the pair $(Y,\cA)$ defined by $\cA_i := \pi_*(\cO(i))$. 

Let $\cB$ be the quasi-coherent sheaf of $\bZ$-graded algebras over $Y$ given by $\cB_i := \pi_* (\cO(i))$. Then we have the following

\bpp  [\cite{Yeu20c}, Proposition 5.32] \label{X_isom_Proj_2}
The sheaf of $\bN$-graded ring $\cB_{\geq 0}$ is Noetherian.
Moreover, there is a canonical isomorphism $\varphi: X \xra{\cong} \Proj^+_Y(\cB)$ over $Y$, together with a canonical isomorphism  $\bigoplus_{i \in \bZ} \varphi^*(\widetilde{\cB(i)}) \xra{\cong} \bigoplus_{i \in \bZ} \cO(i)$ of sheaves of $\bZ$-graded algebras on $X$, such that the composition
\begin{equation*}
\cB_i \xraq{\eqref{M_0_pi_M}} (\pi^+_{\cB})_* (\widetilde{\cB(i)}) \, = \,  \pi_* (\varphi^*(\widetilde{\cB(i)})) \xraq{\cong} \pi_* (\cO(i)) \, =: \, \cB_i
\end{equation*}
is the identity. As a result, $\cB$ is positively $\tfrac{1}{d}$-Cartier.

Furthermore, for any quasi-coherent sheaf $\cF$ on $X$, let $\cM \in \Gr(\cB)$ be the quasi-coherent sheaf of graded $\cB$-module given by $\cM_i := \pi_*( \cF \otimes_{\cO_X} \cO(i))$, then there is a canonical isomorphism $\varphi^*( \widetilde{\cM} ) \xra{\cong} \cF$, such that the composition 
\begin{equation*}
\cM_0 \xraq{\eqref{M_0_pi_M}} (\pi^+_B)_*( \widetilde{\cM}) \, = \,  \pi_* (\varphi^*(\widetilde{\cM})) \xraq{\cong}\pi_*( \cF) \, =: \, \cM_0
\end{equation*}
is the identity.
\epp

Notice that the first part of this proposition says that the setting \eqref{sheaf_of_alg_on_X} always reduces to that of a positively $\tfrac{1}{d}$-Cartier pair $(Y,\cA)$, while the second part answers the question \eqref{tilde_L_relations}(2).

Combining Proposition \ref{M_Ce_L_isoms} and \ref{X_isom_Proj_2}, we have thus answered the questions in \eqref{tilde_L_relations}. This leads to the following version of Serre's equivalence:
\bthm[\cite{Yeu20c}, Theorem 5.34]
\label{Serre_equiv_thm}
If $\cA_{\geq 0}$ is Noetherian and $\tfrac{1}{d}$-Cartier, then the functors $\!\,^0 \! \cL^{+}$ and $(-)^{\sim}$ are inverse equivalences 
\begin{equation}  \label{Serre_equiv_functors}
\begin{tikzcd}
\!\,^0 \! \cL^{+} \,:\, \QCoh(X^+) \ar[r, shift left]
& \QpdGr(\cA) \, : \, (-)^{\sim}  \ar[l, shift left]
\end{tikzcd}
\end{equation}
\ethm

Serre's equivalence gives the relation \eqref{GGM_three_points}(1) we hinted at before. 
Namely, one has the following (see \cite[(3.21)]{Yeu20c})

\bpp  \label{D_IpTR_DQCoh}
If $\cA$ is Noetherian, then there is an inverse pair of exact equivalence
\begin{equation}  \label{QpGrA_IpTR_equiv}
\begin{tikzcd}
\phi^* \, : \, \cD_{\cIpTR}(\Gr(\cA)) \ar[r, shift left]
&  \cD(\QpGr(\cA)) \, : \, {\bm R}\phi_* \ar[l, shift left]
\end{tikzcd}
\end{equation}
\epp

Combining these two results, we see that, under the Cartier assumption, the semi-orthogonal summand $\cD_{\cIpTR}(\Gr(\cA))$ of $\cD(\Gr(\cA))$ is equivalent to $\cD(\QCoh(X^+))$. Without the Cartier assumption, one replaces the projective space $X^+ = \Proj_Y^+(\cA)$ with a stacky version $\PProj_Y^+(\cA)$ (see \cite[Remark 3.19]{Yeu20c}).

The exact equivalence \eqref{QpGrA_IpTR_equiv} restricts to certain bounded coherent subcategories. One has to be careful that the relevant subcategory of $\cD_{\cIpTR}(\Gr(\cA))$ is not given by $\cD_{\cIpTR}(\Gr(\cA)) \cap \Dbcoh(\Gr(\cA))$. Instead, one considers the following
\bdf  \label{IpTR_coh_def}
Assume that $\cA$ is Noetherian. Denote by  $\gr(\cA) \subset \Gr(\cA)$ the full subcategory of sheaves of finitely generated graded modules.
Let $\qpgr(\cA) \subset \QpGr(\cA)$ be the essentially image of $\gr(\cA)$ under $\phi^* : \Gr(\cA) \ra \QpGr(\cA)$.
For each $\spadesuit \in \{\, \, \, ,+,-,b\}$, 
\begin{enumerate}
	\item let $\Dsuit_{\coh}(\QpGr(\cA)) \subset \Dsuit(\QpGr(\cA))$ be the full subcategory consisting of complexes whose cohomology lies in $\qpgr(\cA)$; and
	\item let $\Dsuit_{{\rm coh}(\cIpTR)}(\Gr(\cA)) \subset \Dsuit_{\cIpTR}(\Gr(\cA))$ be the essential image of $\Dsuit_{{\rm coh}}(\Gr(\cA))$ under the functor $\Ce_{\scI^+} : \Dsuit(\Gr(\cA)) \ra \Dsuit_{\cIpTR}(Gr(\cA))$.
\end{enumerate}
\edf

\bpp [\cite{Yeu20c}, Corollary 3.26]  \label{coh_IpTR_QpGr}
For $\spadesuit \in \{-,b\}$, the equivalence \eqref{QpGrA_IpTR_equiv} restricts to give an exact equivalence 
\begin{equation*} 
\begin{tikzcd}
\phi^* \, : \, \Dsuit_{\coh(\cIpTR)}(\Gr(\cA)) \ar[r, shift left]
&  \Dsuit_{\coh}(\QpGr(\cA)) \, : \, {\bm R}\phi_* \ar[l, shift left]
\end{tikzcd}
\end{equation*}
\epp



\section{Log flips and graded rings}  \label{log_flip_subsec}

We continue to work in the setting of \eqref{sheaf_A_setting}, which we will simply refer to as a pair $(Y,\cA)$. The discussion in Section \ref{DGrA_sec} applies for both the positive and the negative directions. This associates two projective morphisms
\begin{equation}  \label{Xm_Y_Xp}
\begin{tikzcd}[row sep = 0]
X^- \, := \, \Proj^-_Y(\cA) \ar[rd, "\pi^-"]  & & \Proj^+_Y(\cA) \, =: \, X^+ \ar[ld, "\pi^+"'] \\
& Y
\end{tikzcd}
\end{equation}

On the scheme $X^- = \Proj^-_Y(\cA)$, we will still write $\cO_{X^-}(i) := \widetilde{A(i)}$. Thus, $\cO_{X^-}(-1)$ is $\pi^-$-ample on $X^-$, perhaps contrary to some conventions.

The main examples of pairs $(Y,\cA)$ that we consider are the pairs associated to log flips between normal varieties. We first recall some standard terminology:

\bdf 
A morphism $f : X \ra Y$ between varieties is said to be a \emph{contraction} if it is projective, and if the canonical map $\cO_Y \ra f_* \cO_X$ is an isomorphism.

A birational contraction $f$ is said to be \emph{small} if the exceptional set $\Ex(f) \subset X$ has codimension $\geq 2$.
\edf

If $X$ and $Y$ are normal varieties, then by Zariski's main theorem, any birational projective morphism $f : X \ra Y$ is a contraction. In this case, for any Weil divisor $D \in \WDiv(X)$, there is a canonical inclusion $f_*\cO(D) \rinto \cO(f_* D)$ if we regard both as subsheaves of $\scK_X = \scK_Y$. One advantage of smallness is that 
\begin{equation}  \label{small_pushforward_div}
\parbox{40em}{If $f : X \ra Y$ is a small birational contraction, then for any Weil divisor $D \in \WDiv(X)$, the canonical inclusion $f_*\cO(D) \rinto \cO(f_* D)$ is an isomorphism.}
\end{equation}


Recall that two birational small contractions
\begin{equation}  \label{Xm_Y_Xp_2}
\begin{tikzcd} [row sep = 0]
X^-  \ar[rd,"\pi^-"] & &  X^+ \ar[ld, "\pi^+"'] \\
& Y
\end{tikzcd}
\end{equation}
of normal varieties over a field $k$ is said to be a \emph{log flip}%
\footnote{The divisor $D^-$ is often written in the form $D^- = K_{X^-} + D^{\prime -}$} if 
there are Weil divisors $D^-$ on $X^-$ and $D^+$ on $X^+$, strict transform of each other, such that
\begin{enumerate}
	\item $-D^-$ is $\bQ$-Cartier and $\pi^-$-ample;
	\item $D^+$ is $\bQ$-Cartier and $\pi^+$-ample.
\end{enumerate}

In this case, if we denote by $D_Y$ their common strict transform to $Y$, then we have the following

\bpp  \label{log_flip_algebra}
The quasi-coherent sheaf of $\bZ$-graded rings $\cA := \bigoplus_{i \in \bZ} \cO_Y(iD_Y)$ is Noetherian. Moreover, the maps \eqref{Xm_Y_Xp} for the resulting pair $(Y,\cA)$ is canonically identified with \eqref{Xm_Y_Xp_2}. Under this identification, there are also canonical identifications of the sheaves of $\bZ$-graded rings 
\begin{equation}   \label{psi_A_tilde_O_iD}
\psi \, : \, 
\bigoplus_{i \in \bZ} \widetilde{\cA(i)}_{X^{\pm}}
\xraq{\cong} 
\bigoplus_{i \in \bZ} \cO_{X^{\pm}}(iD^{\pm})
\end{equation}
Thus, if $dD^-$ is Cartier, then $\cA$ is negatively $\tfrac{1}{d}$-Cartier; if $dD^+$ is Cartier, then $\cA$ is positively $\tfrac{1}{d}$-Cartier. The maps $\psi$ on $X^-$ and $X^+$ are moreover compatible in the sense that the following diagram commutes:
\begin{equation}  \label{O_iD_to_A_tilde_comp_prop}
\begin{tikzcd}
\pi^- _* (\widetilde{\cA(i)}_{X^{-}}) \ar[d, "\pi^-_*(\psi)"', "\cong"]
& \cA_i \ar[l, "\eqref{M_0_pi_M}"'] \ar[r, "\eqref{M_0_pi_M}"] 
\ar[d, equal]
& \pi^+ _* (\widetilde{\cA(i)}_{X^{+}}) 
 \ar[d, "\pi^+_*(\psi)", "\cong"']\\
\pi^- _* \cO_{X^-}(iD^-) 
& \cO_{Y}(iD_Y) \ar[r, equal] \ar[l,equal]
& \pi^- _* \cO_{X^+}(iD^+)
\end{tikzcd}
\end{equation}
where the equalities in the second row are due to \eqref{small_pushforward_div}.
\epp

\bpf
Notice that the conditions in \eqref{sheaf_of_alg_on_X} are satisfied for $\pi^{\pm} : X^{\pm} \ra Y$ and for the sheaves of $\bZ$-graded algebras $\bigoplus_{i \in \bZ} \cO_{X^{\pm}}(iD^{\pm})$ on $X^{\pm}$.
Moreover, by \eqref{small_pushforward_div}, we have canonical isomorphisms $\pi^-_*( \cO(iD^-)) \cong \cA_i \cong \pi^+_*( \cO(iD^+))$, in fact equality as subsheaves of $\scK_Y$. This shows that the sheaf $\cA$ is Noetherian (see Propositions  \ref{Noeth_gr_ring} and \ref{X_isom_Proj_2}).
The claimed isomorphisms $\widetilde{\cA(i)}_{X^{\pm}} \cong \cO_{X^{\pm}}(iD^{\pm})$, the commutativity of \eqref{O_iD_to_A_tilde_comp_prop}, and the claimed $\tfrac{1}{d}$-Cartier property, also follows from Proposition \ref{X_isom_Proj_2}.
\epf

Proposition \ref{log_flip_algebra} shows that for any log flip \eqref{Xm_Y_Xp_2}, the assoicated pair $(Y,\cA)$ determines a log flip in the sense of the following

\bdf  \label{determine_log_flip}
Let $(Y,\cA)$ be a pair as in \eqref{sheaf_A_setting} satisfying $\cA_0 = \cO_Y$. We say that $(Y,\cA)$ \emph{determines a log flip} if the associated diagram \eqref{Xm_Y_Xp} consists of small birational contraction between normal varieties over a field $k$, and if there exists a Weil divisor $D^-$ on $X^-$ and $D^+$ on $X^+$, strict transform of each other, such that there are isomorphisms of sheaves of $\bZ$-graded algebras
\begin{equation}  \label{O_iD_to_A_tilde}
\psi \, : \, 
\bigoplus_{i \in \bZ} \widetilde{\cA(i)}_{X^{\pm}}
\xraq{\cong} 
\bigoplus_{i \in \bZ} \cO_{X^{\pm}}(iD^{\pm})
\end{equation}
which are compatible in the sense that the following diagram commutes:
\begin{equation}  \label{O_iD_to_A_tilde_comp}
\begin{tikzcd}
\pi^- _* (\widetilde{\cA(i)}_{X^{-}}) \ar[d, "\pi^-_*(\psi)"', "\cong"]
& \cA_i \ar[l, "\eqref{M_0_pi_M}"'] \ar[r, "\eqref{M_0_pi_M}"] 
& \pi^+ _* (\widetilde{\cA(i)}_{X^{+}}) 
\ar[d, "\pi^+_*(\psi)", "\cong"']\\
\pi^- _* \cO_{X^-}(iD^-) 
& \cO_{Y}(iD_Y) \ar[r, equal] \ar[l,equal]
& \pi^- _* \cO_{X^+}(iD^+)
\end{tikzcd}
\end{equation}
where the equalities in the second row are due to \eqref{small_pushforward_div}.

Notice that the isomorphism \eqref{O_iD_to_A_tilde} guarantees that $-D^-$ is $\bQ$-Cartier and $\pi^-$-ample, while  $D^+$ is $\bQ$-Cartier and $\pi^+$-ample, so that the associated diagram \eqref{Xm_Y_Xp} is indeed a log flip.

If, moreover, every sheaf $\widetilde{\cA(i)}_{X^{\pm}}$ is maximally Cohen-Macaulay, then we say that the pair $(Y,\cA)$ \emph{determines a Cohen-Macaulay log flip}.
\edf

\brm   \label{A_i_pushforward}
The diagram \eqref{O_iD_to_A_tilde_comp} induces a canonical map $\cA_i \ra \cO_Y(iD_Y)$. In constrast to \eqref{O_iD_to_A_tilde_comp_prop}, we do not require this map to be an isomorphism. In fact, this map is an isomorphism if and only if the sheaf $\cA_i$ is reflexive. Indeed, the commutativity of the diagram \eqref{O_iD_to_A_tilde_comp} with the filled in map $\cA_i \ra \cO_Y(iD_Y)$ shows that this map is an isomorphism on an open subset with complement of codimension $\geq 3$. If $\cA_i$ is reflexive, then this map is an isomorphism by Proposition \ref{refl_sheaf_pushforward} below, since a divisorial sheaf $\cO_Y(iD_Y)$ on a normal variety is always reflexive. The converse is obvious.
\erm

Now we study the inverse problem: given a pair $(Y,\cA)$, deduce some properties for the diagram \eqref{Xm_Y_Xp}. In particular, Proposition \ref{MCM_uple_comp} gives a sufficient condition for $X^{\pm}$ to be Cohen-Macaulay; while Proposition \ref{small_S3_then_log_flip} gives a sufficient condition for $(Y,\cA)$ to determine a log flip.

In the following discussion, we will be mostly interested in the case $\cO_Y = \cA_0$. Otherwise, one may replace $Y$ by $\widetilde{Y} := \Spec_Y \, \cA_0$. 
The first question we are interested in is when $\pi^-$ and $\pi^+$ are birational.
Since this question is local on the base, we may assume that $Y$ is affine.

Thus, let $A$ be a Noetherian $\bZ$-graded ring, and let $R := A_0$. Let $J := (I^- \cdot I^+)_0 \subset R$. It is clear that we have
\begin{equation*}
(I^-)_0 \, =\,  J \, =\,(I^+)_0 
\end{equation*} 

For any graded prime ideal $\mathfrak{p} \subset A$, denote by $\mathfrak{p}_0 \subset R$ its degree zero part. Suppose that $J \not\subset \mathfrak{p}_0$, then there exists some $f \in A_d$, $g \in A_{-d}$, for some $d > 0$, such that $fg \notin \mathfrak{p}_0$. Thus, if we denote by $A_{\mathfrak{p}_0}$ the localization of $A$ with respect to the multiplicative system $S := R \setminus \mathfrak{p}_0$, then we have $(A_{\mathfrak{p}_0})^{(d)} \cong R_{\mathfrak{p}_0}[t,t^{-1}]$.
Thus, the morphisms \eqref{Xm_Y_Xp} become isomorphisms when we localize at $\mathfrak{p}_0 \in \Spec \, R \setminus V(J)$. We summarize this into the first paragraph of the following
\blm  \label{V_J_compl}
For a pair $(Y,\cA)$ as in \eqref{sheaf_A_setting} satisfying $\cA_0 = \cO_Y$, let $\scJ := (\scI^- \cdot \scI^+)_0 \subset \cO_Y$. Then for all $y \in Y \setminus V(\scJ)$, there exists an open neighborhood $y \in U \subset Y \setminus V(\scJ)$ such that $\cA^{(d)}|_U \cong (\cA|_U)_0[t,t^{-1}]$ for some $d>0$, hence the morphisms \eqref{Xm_Y_Xp} are isomorphisms above $Y \setminus V(\scJ)$. 

Therefore, if $\cA$ is furthermore assumed to be a sheaf of integral domains, and if $\cA_{<0} \neq 0$ and $\cA_{>0} \neq 0$, then the morphisms \eqref{Xm_Y_Xp} are birational. 
Moreover, if we let $d := \min\{ 
|i| \in \bZ_{>0} \, | \, i \neq 0 \text{ such that } \cA_i \neq 0 \}$, then there are embeddings of sheaves of $\bZ$-graded rings $\cA \rinto \scK_Y[t,t^{-1}]$ and $\bigoplus_{i \in \bZ} \widetilde{\cA(i)}_{X^{\pm}} \rinto \scK_{X^{\pm}}[t,t^{-1}]$, where $\deg(t) = d$, such that the following diagram commutes:
\begin{equation}  \label{embed_K_comp}
\begin{tikzcd}
\bigoplus_{i \in \bZ} \pi^-_*(\widetilde{\cA(i)}_{X^{-}}) \ar[d, hook]
& \cA \ar[l] \ar[r] \ar[d, hook]
& \bigoplus_{i \in \bZ} \pi^+_*(\widetilde{\cA(i)}_{X^{+}})  \ar[d, hook] \\
\pi^-_* \scK_{X^-}[t,t^{-1}]
& \scK_Y[t,t^{-1}] \ar[l, equal] \ar[r,equal]
& \pi^+_* \scK_{X^+}[t,t^{-1}]
\end{tikzcd}
\end{equation}
\elm

\bpf
The preceding discussion establishes the first paragraph. For the second paragraph, notice that, since $\cA$ is a sheaf of integral domains, and since $\cA_{<0} \neq 0$ and $\cA_{>0} \neq 0$, we have $\cA_i \neq 0$ if and only if $i/d \in \bZ$.
It is then clear that any nonzero homogeneous element $f \in \cA(U)_i$ becomes invertible after passing to the stalk at the generic point $\eta \in Y$. 
In particular, choosing such an element of degree $d$ gives an invertible element $f \in \cA_{\eta} = \cA \otimes_{\cO_Y} \scK_Y$, and hence identifies 
$
\scK_Y[t,t^{-1}]
\cong
\cA \otimes_{\cO_Y} \scK_Y
$
by sending $t$ to $f$.
Since $\cA$ is integral, we have 
$\cA \rinto \cA \otimes_{\cO_X} \scK_X$, which compose to the sought for embedding of $\cA$. 

Notice that $\bigoplus_{i \in \bZ} \widetilde{\cA(i)}_{X^{\pm}}$ is a sheaf of $\bZ$-graded integral domains on $X^{
\pm}$. 
Moreover, recall that there is a canonical map $\cA_i \ra \pi^{\pm}_*(\widetilde{\cA(i)}_{X^{\pm}})$ for each $i \in \bZ$, so that the section $f \in \cA(U)_d$ gives rise to a section $f \in \widetilde{\cA(d)}_{X^{\pm}}((\pi^{\pm})^{-1}(U))$. The same argument then gives the embeddings for $\bigoplus_{i \in \bZ} \widetilde{\cA(i)}_{X^{\pm}}$, making the claimed diagram commutes.
\epf


The preimage of the closed subspace $V(\scJ) \subset Y$ under $\pi^-$ and $\pi^+$ also has a simple description. Indeed, focusing on the affine case, one can show that
$\sqrt{J \cdot A} \, = \, \sqrt{I^- \cdot I^+} $. 
%
%
Localizing to $A_f$ for $\deg(f) > 0$ (resp. $\deg(f) < 0$) then gives 

\blm  \label{V_J_above}
For a pair $(Y,\cA)$ as in \eqref{sheaf_A_setting} satisfying $\cA_0 = \cO_Y$, let $\scJ := (\scI^- \cdot \scI^+)_0 \subset \cO_Y$. Then the preimages of $V(\scJ)$ under the maps \eqref{Xm_Y_Xp} are given by
\begin{equation*}
(\pi^-)^{-1}(V(\scJ)) \, = \, \Proj^-_Y  (\cA/\sqrt{\scI^+})
\qquad \text{and} \qquad 
(\pi^+)^{-1}(V(\scJ)) \, = \, \Proj^+_Y  (\cA/\sqrt{\scI^-})
\end{equation*}
\elm


Now we investigate Cohen-Macaulay property (or more generally the ($S_r$) condition) of $X^-$ and $X^+$.
We start with the following basic result in commutative algebra, the proof of which is left to the reader (see \cite[Tags 02I6, 00GQ, 00GT, 00GU, 0AUK]{Sta}):

\blm  \label{depth_finite_extension}
Let $\varphi : R \ra S$ be an injective finite ring map between Noetherian rings. If $S$ is catenary of finite Krull dimension then so is $R$, and we have $\dim(R_{\mathfrak{p}}) = \dim(S_{\mathfrak{q}})$ for all prime ideal $\mathfrak{q} \in \Spec \, S$ above $\mathfrak{p} \in \Spec \, R$.
Under this catenary condition,  the following two conditions are equivalent:
\begin{enumerate}
	\item $S$ satisfies Serre's condition $(S_r)$, \ie, ${\rm depth}_{S_{\mathfrak{q}}}(S_{\mathfrak{q}}) \geq \min\{ r , \dim(S_{\mathfrak{q}}) \}$ for all prime ideals $\mathfrak{q} \in \Spec \, S$.
	\item $S$ is ``maximally $(S_r)$" as a module over $R$, \ie, it satisfies the condition ${\rm depth}_{R_{\mathfrak{p}}}(S_{\mathfrak{p}}) \geq \min\{ r , \dim(R_{\mathfrak{p}}) \}$ for all prime ideals $\mathfrak{p} \in \Spec \, R$, where we denote by $S_{\mathfrak{p}}$ the localization of $S$ with respect to the multiplicative subset $R \setminus \mathfrak{p}$
\end{enumerate}
\elm


\bcor  \label{S_r_cond_Q_Car}
Let $B$ be a Noetherian $\bZ$-graded ring such that $B^{(d)} \cong B_0[t,t^{-1}]$ for some $d > 0$. 
Suppose that $B_0$ is universally catenary. Then the following two conditions are equivalent:
\begin{enumerate}
	\item $B$ satisfies Serre's condition $(S_r)$.
	\item Each $B_i$, $i \in \bZ$ is ``maximally $(S_{r-1})$" as a module over $B_0$, \ie, it satisfies the condition ${\rm depth}_{(B_0)_{\mathfrak{p}}}((B_i)_{\mathfrak{p}}) \geq \min\{ r-1 , \dim((B_0)_{\mathfrak{p}}) \}$ for all prime ideals $\mathfrak{p} \in \Spec \, B_0$.
\end{enumerate}
In particular, $B$ is Cohen-Macaulay of dimension $n+1$ if and only if $B_0$ is Cohen-Macaulay of dimension $n$ and each $B_i$ is maximal Cohen-Macaulay as a module over $B_0$.
\ecor

\bpf
Apply Lemma \ref{depth_finite_extension} to the finite extension $B^{(d)} \rinto B$, which is valid because $B_0$ is assumed to be universally catenary. 
By assumption, we have $B^{(d)} \cong B_0[t,t^{-1}]$. Moreover, by considering the $\bZ$-grading mod $d$, we see that, as a module over $B^{(d)}$, $B$ splits as a direct sum $B = \bigoplus_{i = 0}^{d-1} M_{(i)}$, where $M_{(i)} := \bigoplus_{j \in \bZ} \, B_{i + jd} \cong B_i [t,t^{-1}]$. 
Since the depth of a finite direct sum is equal to the minimal of the depths of the summands, we see that in this case, condition (2) in Lemma \ref{depth_finite_extension} says that each $M_{(i)}$ is maximally $(S_r)$ over $B^{(d)}$.
As we have observed above, we have $(B^{(d)},M_{(i)}) \cong (B_0[t,t^{-1}], B_i[t,t^{-1}])$, so that this last condition is furthermore equivalent to the condition that each $B_i$ is maximally $(S_{r-1})$ over $B_0$. Thus, the result follows from Lemma \ref{depth_finite_extension}.
\epf

This can be used to prove the following

\bpp  \label{MCM_uple_comp}
If $A$ is a Cohen-Macaulay $\bZ$-graded Noetherian ring of (pure) dimension $n+1$, where $A_0$ is universally catenary, then both $X^-$ and $X^+$ are Cohen-Macaulay of (pure) dimension $n$. Moreover, for each $i \in \bZ$, the sheaf $\cO_{X^{\pm}}(i) := \widetilde{A(i)}_{X^{\pm}}$ is maximally Cohen-Macaulay. 
\epp

\bpf
It suffices to prove this for $X^+$. For any homogeneous element $f \in A_d$ of degree $d>0$, the localization $B := A_f$ is still Cohen-Macaulay of pure dimension $n+1$, and $B_0$ is still universally catenary. Thus, we may apply Corollary \ref{S_r_cond_Q_Car} to obtain the desired result.
\epf

%


We also apply Corollary \ref{S_r_cond_Q_Car} to give sufficient conditions for a $\bZ$-graded ring to give rise to a log flip.
We start with the following
\blm  \label{A_S3_then_refl}
Suppose that the $\bZ$-graded ring $A$ is a normal domain, then the schemes $X^-$, $X^+$ and $Y := \Spec \, A_0$ are normal integral schemes. If $A$ moreover satisfies Serre's condition $(S_3)$, and if $A_{<0} \neq 0$ and $A_{>0} \neq 0$, then for each $i \in \bZ$, the sheaf $\cO_{X^{\pm}}(i) := \widetilde{A(i)}_{X^{\pm}}$ is reflexive. 
\elm

\bpf
Recall that a domain is normal if and only if it is integrally closed. From this, one sees easily that
\begin{equation}  \label{normal_domain_deg_zero}
\parbox{40em}{If $A$ is a $\bZ$-graded normal domain, then $A_0$ is also a normal domain.}
\end{equation}

Since $X^+$ is locally given by $\Spec \, (A_{f})_0$, and since localization preserves normality, the first statement for $X^+$ follows directly from \eqref{normal_domain_deg_zero}. Similarly for $X^-$. The statement for $Y$ also follows from \eqref{normal_domain_deg_zero}. 

If $A$ satisfies Serre's condition $(S_3)$, then so does $A_f$, for $\deg(f) > 0$. Applying Corollary \ref{S_r_cond_Q_Car} to $B := A_f$, we see for each $i \in \bZ$, the module $(A_f)_i$ over $(A_f)_0$ is maximally $(S_2)$. In other words, each of the sheaves $\widetilde{\cA(i)}_{X^+}$ is maximally $(S_2)$. Since $X^+$ is normal, this implies that $\widetilde{\cA(i)}_{X^+}$ is reflexive (see, e.g., \cite[Proposition 1.4.1(b)]{BH93}). 
\epf

\bpp  \label{small_S3_then_log_flip}
Let $(Y,\cA)$ be a pair as in \eqref{sheaf_A_setting}. Suppose that
\begin{enumerate}
	\item $Y$ is a variety over a field $k$, and $\cA_0 = \cO_Y$.
	\item $\cA$ is a sheaf of normal domains satisfying the condition $(S_3)$.
	\item The closed subsets $V(\scI^-)$ and $V(\scI^+)$ of $\Spec_Y \, \cA$ both have codimension $\geq 2$.
\end{enumerate}
then $(Y,\cA)$ determines a log flip in the sense of Definition \ref{determine_log_flip}. If the sheaf $\cA$ of algebras is furthermore Cohen-Macaulay, then it determines a Cohen-Macaulay log flip.
\epp

\bpf
Normality of $X^-$, $X^+$ and $Y$ follows from Lemma \ref{A_S3_then_refl}. 
The assumption (3) implies that $\cA_{<0} \neq 0$ and $\cA_{>0} \neq 0$, so that Lemma \ref{V_J_compl} shows in particular that $\pi^-$ and $\pi^+$ are birational, and $\dim(Y) = \dim(\cA(U)) - 1$ for any affine open $U \subset Y$. 
Since $\dim(\Proj(B)) \leq \dim(B) - 1$ for any Noetherian $\bN$-graded ring $B$, we see from Lemmae \ref{V_J_compl} and \ref{V_J_above} that the exceptional loci of $\pi^-$ and $\pi^+$ have codimension $\geq 2$, proving the first condition of Definition \ref{determine_log_flip}.

For the second condition, apply Lemma \ref{V_J_compl} again, which gives embeddings $\bigoplus_{i \in \bZ} \widetilde{\cA(i)}_{X^{\pm}} \rinto \scK_{X^{\pm}}[t,t^{-1}]$ of $\bZ$-graded $k$-algebras. Since $\widetilde{\cA(i)}_{X^{\pm}}$ are reflexive sheaves (see Lemma \ref{A_S3_then_refl}), these embeddings identify $\widetilde{\cA(i)}_{X^{\pm}} = \cO(D_i^{\pm})$, for some Weil divisor $D_i^{\pm}$ on $X^{\pm}$. Since the embedding into  $\scK_{X^{\pm}}[t,t^{-1}]$ is multiplicative, we have $D_i^{\pm} + D_j^{\pm} \leq D^{\pm}_{i+j}$ for all $i,j \in \bZ$. On the other hand, if we choose $d > 0$ such that the pair $(Y,\cA)$ is $\tfrac{1}{d}$-Cartier, then we have $D_i^{\pm} + D_d^{\pm} = D^{\pm}_{i+d}$. Thus, we must have $D_i^{\pm} + D_j^{\pm} = D^{\pm}_{i+j}$ for all $i,j \in \bZ$. In particular, we have $D_i^{\pm} = iD^{\pm}$ for $D^{\pm} := D_1^{\pm}$.
This gives the desired isomorphisms \eqref{O_iD_to_A_tilde}. 
The commutativity of \eqref{O_iD_to_A_tilde_comp} then follows from that of \eqref{embed_K_comp}. Therefore, it suffices to show that $D^+$ is the strict transform of $D^-$. 
Indeed, the subsheaves $\cO_{X^-}(D^-)$ and $\cO_{X^+}(D^+)$ of $\scK_{X^-}$ and $\scK_{X^+}$ are defined as the images of the embedding $\widetilde{\cA(i)}_{X^{\pm}} \rinto \scK_{X^{\pm}}$. By the commutativity of \eqref{embed_K_comp}, they coincide on the canonically identified open subschemes
\begin{equation*}
X^- \setminus \Proj^-_Y  (\cA/\sqrt{\scI^+})  \xraq[\cong]{\pi^-}  
Y \setminus V(\scJ) 
\xlaq[\cong]{\pi^+}  X^+ \setminus \Proj^+_Y  (\cA/\sqrt{\scI^-}) 
\end{equation*}
(see Lemmae \ref{V_J_compl} and \ref{V_J_above}). Since these open subschemes have complement of codimension $\geq 2$, the Weil divisors $D^-$ and $D^+$ are strict transforms of each other.

For the last statement, simply apply Proposition \ref{MCM_uple_comp}.
\epf

\brm  \label{small_then_log_flip_remark}
Without the $(S_3)$ assumption on $\cA$ in Proposition \ref{small_S3_then_log_flip}(2), one can also show that the $d$-uple component $(Y,\cA^{(d)})$ determines a log flip for any $d > 0$ such that $(Y,\cA)$ is $\tfrac{1}{d}$-Cartier. This is, in effect, the proof of \cite[Proposition 1.6]{Tha96}. However, it seems that in order to establish \cite[Proposition 1.6]{Tha96} as is stated there, some condition similar our $(S_3)$ assumption is necessary.
\erm

\section{Flips and flops}  \label{flip_flop_subsec}

%



In this section, we will specify some properties of Grothendieck duality associated to flips and strongly crepant flops (see Definition \ref{flip_flop_def}) between Cohen-Macaulay normal varieties.
We start by recalling some basic notions about the functor $f^!$, following \cite[Tag 0DWE]{Sta}. Recall that we work with Conventions \ref{Dcat_conv1} and \ref{Dcat_conv2} throughout.


Recall from \cite[Tag 0F42]{Sta} that, given a morphism $f : X \ra Y$ of schemes, both separated and of finite type over a separated Noetherian base scheme $S$, then we may define the upper shriek functor%
\footnote{See Convention \ref{Dcat_conv2}.}
$f^! : \cD(\QCoh(Y)) \ra \cD(\QCoh(X))$ as follows: first factor $f$ as the composition 
$X \xra{j} \bar{X} \xra{\tilde{f}} Y$ of an open immersion $j$ and a proper morphism $\tilde{f}$, both over $S$, then define $f^!$ to be the composition $f^! = j^* \circ \tilde{f}^!$, where $\tilde{f}^!$ is right adjoint to ${\bm R}\tilde{f}_*$, and $j^*$ is the restriction.
A different factoring $f = j' \circ \tilde{f}'$ gives rise to canonically isomorphic functors $j^{\prime *} \circ \tilde{f}^{\prime !} \cong j^* \circ \tilde{f}^!$.

We also recall (see \cite[Tag 0B6N]{Sta}) that if $f : X \ra Y$ is proper, then there is a canonical map 
\begin{equation}  \label{shriek_tensor}
f^!(\cF) \otimes_{\cO_X}^{{\bm L}} {\bm L}f^* \cG \raq f^!( \cF \otimes_{\cO_Y}^{{\bm L}} \cG  )
\end{equation}
in $\cD(\QCoh(X))$ for each $\cF,\cG \in \cD(\QCoh(Y))$, which is adjoint to the map
\begin{equation*}
{\bm R}f_*( f^!(\cF) \otimes_{\cO_X}^{{\bm L}} {\bm L}f^* \cG ) \, = \, {\bm R}f_* f^!(\cF) \otimes_{\cO_Y}^{{\bm L}} \cG \xraq{\epsilon \otimes \id} \cF \otimes_{\cO_Y}^{{\bm L}} \cG
\end{equation*}
where $\epsilon$ is the adjunction counit.
Moreover, \eqref{shriek_tensor} is an isomorphism if $\cG$ is perfect.

A crucial feature of the upper shriek functor is that it sends dualizing complexes to dualizing complexes (see \cite[Tags 0A7B, 0A87, 0AA3]{Sta}).
In particular, for a scheme $X$ separated and of finite type over a field $k$, if we denote by $\pi : X \ra \Spec \,k$ the projection map, then $\omega_X^{\bullet} := \pi^!(k)$ is a dualizing complex, called the \emph{canonical dualizing complex} of $X$. 
The following Lemma gives the correct shift in order to normalize the restriction $(\omega_X^{\bullet})_x$ to each stalk:
\blm
For any point $x\in X$, if we let $\Delta_x \in \bZ$ be the unique integer such that $(\omega_X^{\bullet})_x[-\Delta_x]$ is a normalized dualizing complex over $\cO_{X,x}$ in the sense of \cite[Tag 0A7M]{Sta}, then we have $\Delta_x = \dim(\overline{\{x\}})$.
\elm

\bpf
If $x$ is a closed point then this follows from the adjunctions
\begin{equation*}
\begin{split}
& \RHom_{\cO_{X,x}}(k(x), (\omega_X^{\bullet})_x) \, \simeq \, 
\ \bigl( \, \RcHom_{\cO_X}((i_x)_* k(x), \omega_X^{\bullet}) \, \bigr)_x \\ 
& \quad \simeq \, 
 \bigl( \, (i_x)_*  \RHom_{k(x)} (k(x), i_x^! \omega_X^{\bullet}) \, \bigr)_x
 \, \simeq \, \RHom_k(k(x),k)
\end{split}
\end{equation*}
where $i_x$ is the closed immersion $i_x : \Spec \, k(x) \ra X$, and the last quasi-isomorphism follows from the functoriality $\pi_x^! \cong i_x^! \pi^!$
for the composition $\pi_x : \Spec \, k(x) \xra{i_x} X \xra{\pi} \Spec \, k$. 
The general case then follows from the fact that $x\mapsto \Delta_x$ is a dimension function on $X$ (see \cite[Tag 0A7Z]{Sta}).
\epf
%
As a result, for any coherent sheaf $\cF \in \Coh(X)$ of $\dim({\rm supp}(\cF_x)) = d_x(\cF)$ and ${\rm depth}(\cF_x) = \delta_x(\cF)$, we have by \cite[Tag 0A7U]{Sta}
\begin{enumerate}
	\item If $\cH^i( \RcHom_{\cO_X}(\cF, \omega_X^{\bullet}))_x \neq 0$ then $-\Delta_x-d_x(\cF) \leq i \leq -\Delta_x-\delta_x(\cF)$;
	\item $\cH^{-\Delta_x - \delta_x(\cF)}( \RcHom_{\cO_X}(\cF, \omega_X^{\bullet}))_x \neq 0$.
\end{enumerate}
In particular, if $X$ has pure dimension $n$ (and is catenary since it is of finite type over a field), then these numbers may be rewritten as
\begin{equation*}
-\Delta_x-d_x(\cF) = -n + (d_x(\cO_x) - d_x(\cF))
\qquad \text{and} \qquad 
-\Delta_x-\delta_x(\cF) = -n + (d_x(\cO_x) - \delta_x(\cF))
\end{equation*}
so that in particular $\omega_X^{\bullet}$ itself has cohomology concentrated in degrees $\geq -n$.
Moreover, $\cF$ is maximal Cohen-Macaulay if and only if $\RcHom_{\cO_X}(\cF, \omega_X^{\bullet})$ has cohomology concentrated in degree $-n$.

If $X$ is moreover proper over $\Spec \, k$, then $\pi^!$ is right adjoint to ${\bm R}\pi_*$, so that for each $\cF \in \Dbcoh(X)$, there is a canonical isomorphism
\begin{equation}  \label{Groth_duality_k_1}
(\bH^i(X , \cF))^* \, \cong \, \EExt^{-i} (\cF , \omega_X^{\bullet} ) \, = \, 
\bH^{-i}(X , \RcHom_{\cO_X}(\cF, \omega_X^{\bullet}))
\end{equation}
Thus, in particular, if we let $\omega_X := \cH^{-n}( \omega_X^{\bullet} )$, then for any coherent sheaf $\cF \in \Coh(X)$, the lowest nonzero degree for which $\bH^{\bullet}(X , \RcHom(\cF, \omega_X^{\bullet}))$ could possibly be nonzero happens at $\bullet = -n$, where it is given by
$H^{0}(X , \cHom(\cF, \omega_X)) = \Hom_X(\cF,\omega_X)$.
At this degree, \eqref{Groth_duality_k_1} becomes
\begin{equation*}
H^n(X,\cF)^* \, \cong \, \Hom_{\cO_X}(\cF,\omega_X)
\end{equation*}
Thus, when $X$ is proper over $k$, $\omega_X := \cH^{-n}( \omega_X^{\bullet} )$ represents the functor $\Coh(X)^{\op} \ra \Mod(k)$, $\cF \mapsto H^n(X,\cF)^*$, and is therefore unique up to canonical isomorphism, known as the \emph{dualizing sheaf on $X$}.

If $X$ is not proper over $k$, the coherent sheaf $\omega_X := \cH^{-n}( \omega_X^{\bullet})$ is still well-defined, although it may not play the role of a dualizing sheaf. 
On any normal variety $X$, it is in fact isomorphic to $\cO_X(K_X)$ (see Theorem \ref{O_KX_dualizing} below). This follows by restricting to $X$ the corresponding isomorphism \cite[Proposition 5.75]{KM98} on a compactification $X \subset \overline{X}$. However, since we will need some details in the construction of such an isomorphism, we provide some details to the second proof hinted at in \cite{KM98}, written in a form more cogenial to derived categories. We start by recalling some preparatory results:

\bpp  \label{omega_depth}
If $X$ is a scheme of pure dimension $n$, separated and of finite type over $k$, then the stalks of the coherent sheaf $\cH^{-n}( \omega_X^{\bullet})$ satisfy ${\rm depth}_{\cO_{X,x}}(\cH^{-n}( \omega_X^{\bullet})_x) \geq \min\{2, \dim(\cO_{X,x})\}$.
\epp

\bpf
We have seen above that the stalk $\cH^{-n}( \omega_X^{\bullet})_x$ is again the lowest degree cohomology of the normalized dualizing complex $(\omega_X^{\bullet})_x[-\Delta_x]$ over $\cO_{X,x}$.
Thus the statement follows from (the proof of) \cite[Tag 0AWE]{Sta}. 
\epf

\bcor  \label{omega_reflexive}
If $X$ is a scheme of pure dimension $n$, separated and of finite type over $k$, and if $\cO_{X,x}$ is Gorenstein for all $x\in X$ with ${\rm depth}(\cO_{X,x}) \leq 1$, then $\cH^{-n}( \omega_X^{\bullet})$ is reflexive. For example, this second condition is satisfied if $X$ is normal.
\ecor

\bpf
This follows from Proposition \ref{omega_depth} and the usual characterization of reflexive coherent sheaves by depth condition (see, e.g., \cite[Proposition 1.4.1(b)]{BH93}).
\epf

\bpp[\cite{Har80}, Proposition 1.6]  \label{refl_sheaf_pushforward}
If $\cF$ is a reflexive coherent sheaf on a normal integral scheme $X$, then for any closed subset $Z \subset X$ of codimension $\geq 2$, the canonical map $\cF \ra j_*j^* \cF$ is an isomorphism, where $j : X \setminus Z \rinto X$ is the inclusion map.
\epp

Now, suppose we are given any proper morphism $f : W \ra X$ between schemes of pure dimension $n$, separated and of finite type over $k$, then the adjunction $({\bm R}f_* \vdash f^!)$, together with the canonical identification $\omega_W^{\bullet} = f^! \omega_X^{\bullet}$, determines a counit morphism 
\begin{equation}   \label{Tr_W_X}
\Tr_{W/X} \, : \, {\bm R}f_* (\omega_W^{\bullet}) \raq \omega_X^{\bullet}
\end{equation}
so that if we take the lowest degree cohomology sheaves, we have a map of coherent sheaves
\begin{equation}  \label{Tr0_W_X}
\!\,^0 \Tr_{W/X} \, : \, f_* \cH^{-n}(\omega_W^{\bullet}) \, \cong \, 
\cH^{-n}( {\bm R}f_* (\omega_W^{\bullet}) )  \xraq{\cH^{-n}(\Tr_{W/X})} \cH^{-n}(\omega_X^{\bullet})
\end{equation}

This can be used to prove the following

\bthm  \label{O_KX_dualizing}
Assume ${\rm char}(k) = 0$, then on any normal variety $X$ over $k$ of dimension $n$, there is an isomorphism $\Phi_X : \cO(K_X) \xra{\cong} \cH^{-n}( \omega_X^{\bullet} ) $.

\ethm

\bpf
Choose a resolution of singularities $f : W \ra X$. 
Since $W$ is smooth, there is an isomorphism $\Phi_W : \cO(K_W) \xra{\cong} \cH^{-n}( \omega_W^{\bullet} )$. 
Then consider the diagram
\begin{equation} \label{Phi_Tr_comm_diag} 
\begin{tikzcd}
f_* \cO(K_W) \ar[r, "f_*(\Phi_W)"] \ar[d, hook]  &  f_* \cH^{-n}( \omega_W^{\bullet} ) \ar[d, "\!\,^0 \Tr_{W/X} "] \\
\cO(K_X) \ar[r, dashed, "\Phi_X"] &  \cH^{-n}( \omega_X^{\bullet} )
\end{tikzcd}
\end{equation}

Let $Z = f(\Ex(f))$, which has codimension $\geq 2$ since $f$ is a birational contraction. 
Since $f$ is an isomorphism outside $Z$, both of the vertical maps are isomorphisms over $X\setminus Z$, so that the dashed arrow in \eqref{Phi_Tr_comm_diag} exists uniquely on $X\setminus Z$. 
Since both $\cO(K_X)$ and $\cH^{-n}( \omega_X^{\bullet} )$ are reflexive (see Corollary \ref{omega_reflexive}), we see by Proposition \ref{refl_sheaf_pushforward} that the dashed arrow $\Phi_X$ in \eqref{Phi_Tr_comm_diag} exists uniquely over $X$.
\epf

We also recall the notion of rational singularities:
\bdf
Let $X$ be a variety over a field of characteristic zero. Then $X$ is said to have \emph{rational singularities} if for some, and hence any, resolution of singularities $f : W \ra X$, the canonical map
\begin{equation}  \label{OX_OW}
\cO_X \raq {\bm R}f_*(\cO_W)
\end{equation}
is an isomorphism.
This is equivalent to the condition that $X$ is normal and ${\bm R}^if_*(\cO_W) = 0$ for all $i > 0$.
\edf

Notice that the maps \eqref{Tr_W_X} and \eqref{OX_OW} are dual to each other. Namely, if we apply $\RcHom_{\cO_X}(-,\omega_X^{\bullet})$ to \eqref{OX_OW}, then we have a map
\begin{equation*}
{\bm R}f_* (\omega_W^{\bullet}) 
\, \cong \,
\RcHom_{\cO_X}({\bm R}f_*\cO_W ,\omega_X^{\bullet})
\raq \RcHom_{\cO_X}(\cO_X,\omega_X^{\bullet}) = \omega_X^{\bullet}
\end{equation*}
where the first isomorphism is the local adjunction isomorphism. This map can be easily shown to coincide with \eqref{Tr_W_X}. 
%
This can be used to give a simple proof of the following result (see also \cite[Theorem 5.10]{KM98}):
\bpp  \label{rat_sing_equiv}
Let $X$ be a variety over a field of characteristic zero, and let $f : W \ra X$ be a resolution of singularities. Then the followings are equivalent:
\begin{enumerate}
	\item $X$ has rational singularities.
	\item The map \eqref{Tr_W_X} is an isomorphism.
	\item $X$ is Cohen-Macaulay and the map \eqref{Tr0_W_X} is an isomorphism.
\end{enumerate}
\epp

\bpf
The equivalence between (1) and (2) follows from the fact that \eqref{Tr_W_X} and \eqref{OX_OW} are dual to each other.
For the equivalence between (2) and (3), recall the Grauert-Riemenschneider vanishing theorem \cite{GR70} (see also \cite[Corollary 2.68]{KM98}, or the proof of Corollary \ref{KV_vanishing_2} below):
\begin{equation}  \label{Grauert-Riemenschneider}
\parbox{40em}{If $f : W \ra X$ is birational morphism between projective varieties over a field of characteristic zero, and $W$ is smooth, then ${\bm R}^i f_* (\omega_W) = 0$ for all $i > 0$}
\end{equation}
Thus, the map \eqref{Tr_W_X} is an isomorphism if and only if $\omega_X^{\bullet}$ has cohomology concentrated in degree $-n$ (\ie, $X$ is Cohen-Macaulay), and the map \eqref{Tr0_W_X} is an isomorphism.
\epf

The condition of rational singularities allows one to extend some results over smooth varieties to varieties with rational singularities. An example is a relative form of the Kawamata-Viehweg vanishing theorem (see Corollary \ref{KV_vanishing_2} below). We first recall the following (see, \eg, \cite[Theorem 2.64]{KM98})

\bthm  \label{KV_vanishing_1}
If $X$ is a smooth and projective variety over a field of characteristic zero, and if $\scL \in \Pic(X)$ is big and nef, then we have $H^i(X, \omega_X \otimes \scL) = 0$ for all $i > 0$.
\ethm

\bcor  \label{KV_vanishing_2}
Let $\pi : X \ra Y$ be a proper morphism between varieties over a field of characteristic zero, where $X$ is projective with at most rational singularities. Let $\scL \in \Pic(X)$. Assume either
\begin{enumerate}
	\item $\pi$ is a birational contraction, and $\scL$ is nef; or
	\item $\scL$ is $\pi$-ample
\end{enumerate}
then we have ${\bm R}^i \pi_*(\omega_X \otimes \scL) = 0$ for all $i > 0$.
\ecor

\bpf
Choose any ample line bundle $\scH$ on $Y$. Recall that a Leray spectral sequence shows that
\begin{equation}  \label{vanishing_ample_enough}
\parbox{40em}{If $\pi : X \ra Y$ is a proper morphism, then for any $\cF \in \Coh(X)$, we have ${\bm R}^i\pi_*(\cF) = 0$ for all $i > 0$ if and only if there exists $n_0 \geq 0$ such that $H^i(X , \cF \otimes \pi^* \scH^{\otimes n}) = 0$ for all $n \geq n_0$.}
\end{equation}

Notice that in either case of our assumption, there exists $n_0 \geq 0$ such that $\scL \otimes \pi^*(\scH^{\otimes n})$ is big and nef for all $n \geq n_0$. 
By the isomorphism ${\bm R}^i \pi_*(\omega_X \otimes \scL \otimes \pi^*(\scH^{\otimes n})) \cong {\bm R}^i\pi_*(\omega_X \otimes \scL) \otimes \scH^{\otimes n}$, the result therefore follows from Theorem \ref{KV_vanishing_1} if $X$ is assumed to be smooth. In particular, this justifies \eqref{Grauert-Riemenschneider}.

To extend this to the case when $X$ has rational singularities, choose a resolution of singularities $f : W \ra X$, and let $g = \pi \circ f$. Then $f^*(\scL \otimes \pi^*(\scH^{\otimes n})) = f^*\scL \otimes g^*(\scH^{\otimes n})$ remains big and nef for all $n \geq n_0$. Applying \eqref{vanishing_ample_enough} to $g : W \ra Y$, it therefore follows from Theorem \ref{KV_vanishing_1} that ${\bm R}^i g_* ( \omega_W \otimes f^*\scL ) = 0$ for all $i > 0$.
We have seen in Proposition \ref{rat_sing_equiv} that $X$ having rational singularities implies that ${\bm R} f_* ( \omega_W) \cong \omega_X$, and hence ${\bm R} f_* ( \omega_W \otimes f^*\scL ) \cong \omega_X \otimes \scL$. Thus, the condition ${\bm R}^i g_* ( \omega_W \otimes f^*\scL ) = 0$ can be rewritten as ${\bm R}^i \pi_* ( \omega_X \otimes \scL ) = 0$.
\epf

Assume from now on ${\rm char}(k) = 0$. Then Theorem \ref{O_KX_dualizing} allows us to formulate a crucial homological property of some classes of flips and flops.
Suppose we are given a diagram of birational contractions between normal varieties, as in \eqref{Xm_Y_Xp_2},
which satisfies the condition
\begin{equation}  \label{pseudo_can_cond}
\parbox{40em}{The canonical inclusion maps $\pi^-_*(\cO(K_{X^-})) \rinto \cO(K_Y)$ and $\pi^+_*(\cO(K_{X^+})) \rinto \cO(K_Y)$ are isomorphisms.}
\end{equation}
This is automatic if both $\pi^-$ and $\pi^+$ are small, in view of \eqref{small_pushforward_div}. It is also true if $X^-$ and $X^+$ are quasi-Calabi-Yau, \ie, if $K_{X^-} \equiv 0$ and $K_{X^+} \equiv 0$.

Using \eqref{pseudo_can_cond}, we have a canonical map
\begin{equation}  \label{pseudo_can_map}
\cO(K_Y) \, = \, \pi^+_* \cO(K_{X^+}) \raq  {\bm R}\pi^+_* \cO(K_{X^+})
\end{equation}

Now, find a common resolution of singularities
\begin{equation}  \label{common_reosln_W}
\begin{tikzcd} [row sep = 0]
& W \ar[ld, "f^-"'] \ar[rd, "f^+"] & \\
X^- \ar[rd, "\pi^-"] & & X^+ \ar[ld, "\pi^+"'] \\
& Y &
\end{tikzcd}
\end{equation}
and use $f^+$ and $f^-$ to define the maps $\Phi_{X^-}$, $\Phi_{X^+}$ and $\Phi_{Y}$
as in the proof of Theorem \ref{O_KX_dualizing}. Since these are defined as the unique map making the diagram \eqref{Phi_Tr_comm_diag} commute, we see easily that the following diagram commutes:
\begin{equation}  \label{Phi_maps_commute}
\begin{tikzcd}
\pi^-_* \cO(K_{X^-})  \ar[d, "\pi^-_*(\Phi_{X^-})"'] \ar[rr, equal]
& & \cO(K_{Y})  \ar[d, "\Phi_{Y}"]
& & \pi^+_* \cO(K_{X^+}) \ar[d, "\pi^+_*(\Phi_{X^+})"] \ar[ll, equal]  \\
\pi^-_* \cH^{-n}(\omega_{X^-}^{\bullet}) \ar[rr, "\!\,^0 \Tr_{X^-/Y}"]
& &  \cH^{-n}(\omega_{Y}^{\bullet}) 
& & \pi^+_* \cH^{-n}(\omega_{X^+}^{\bullet})  \ar[ll, "\!\,^0 \Tr_{X^+/Y}"']
\end{tikzcd}
\end{equation}

Now suppose that both $X^-$ and $X^+$ are Cohen-Macaulay, then there is an isomorphism
\begin{equation}  \label{pi_shriek_KX}
\cO(K_{X^+}) \xra[\cong]{\Phi_{X^+}} \cH^{-n}(\omega_{X^+}^{\bullet}) \xra[\cong]{\iota} \omega_{X^+}^{\bullet}[-n] = (\pi^+)^! \omega_{Y}^{\bullet}[-n]
\qquad \text{ in } \Dbcoh(\QCoh(X^+)) 
\end{equation}
and similarly for $X^-$. Likewise, there is also a map
\begin{equation}  \label{O_K_Y_omega_map}
\cO(K_Y) \xraq[\cong]{\Phi_Y} \cH^{-n}( \omega_{Y}^{\bullet}[-n] ) \xraq{\iota} \omega_{Y}^{\bullet}[-n]
\qquad \text{ in } \Dbcoh(\QCoh(Y)) 
\end{equation}
which is an isomorphism if $Y$ is Cohen-Macaulay.

We may use the isomorphism \eqref{pi_shriek_KX} to express Grothendieck duality for the morphism $\pi^+$ in terms of the sheaf $\cO(K_{X^+})$ and the dualizing complex $\omega_{Y}^{\bullet}[-n]$. 
Namely, for any $\cF \in \Dmcoh(\QCoh(X^+))$, there is a local adjunction isomorphism%
\footnote{Recall that we work under Conventions \ref{Dcat_conv1} and \ref{Dcat_conv2}. Thus, we impose the condition $\cF \in \Dmcoh(\QCoh(X^+))$ to guarantee that our local adjunction isomorphisms coincide with the usual ones on $\RcHom^{\clubsuit}(-,-)$.}
\begin{equation}  \label{local_adj_map_Phi}
{\bm R}\pi^+_* \RcHom_{\cO_{X^+}}(\cF , \cO(K_{X^+})) \xraq{\cong} 
\RcHom_{\cO_Y}( {\bm R}\pi^+_* \cF , \omega_{Y}^{\bullet}[-n] )
\end{equation}
which sends an element $\varphi \in {\bm R}\Gamma(U ; {\bm R}\pi^+_* \RcHom_{\cO_{X^+}}(\cF , \cO(K_{X^+})) ) \simeq \RHom_{(\pi^+)^{-1}U}( \cF , \cO(K_{X^+}) )$ to the composition (defined on $U$)
\begin{equation}  \label{local_adj_map_Phi_descr}
{\bm R}\pi^+_*\cF \xra{{\bm R}\pi^+_*(\varphi)} {\bm R}\pi^+_*(\cO(K_{X^+})) 
\xra{ {\bm R}\pi^+_*( \ref{pi_shriek_KX} ) } {\bm R}\pi^+_*(\omega_{X^+}^{\bullet}[-n])
\xra{ \Tr_{X^+/Y} } \omega_{Y}^{\bullet}[-n]
\end{equation}


The commutativity of  \eqref{Phi_maps_commute} can be used to establish a compatibility condition between the local adjunction isomorphism \eqref{local_adj_map_Phi} and the corresponding version for $\pi^-$.
Namely, the commutativity of the right square of \eqref{Phi_maps_commute} may be reformulated as the following Lemma, which will be used in the next section in the proof of Theorem \ref{HFF_main_examples}.
\blm  \label{compatible_lemma_1}
The composition
\begin{equation}  \label{compatible_lemma_1_map_1}
\begin{split}
\cO_Y &\raq {\bm R}\pi^+_* \cO_{X^+} \xraq{ {\bm R}\pi^+_*(\Delta) }
{\bm R}\pi^+_* \RcHom_{\cO_{X^+}}(\cO(K_{X^+}) , \cO(K_{X^+})) \\
&\xraq{ \eqref{local_adj_map_Phi} } 
\RcHom_{\cO_Y}( {\bm R}\pi^+_* \cO(K_{X^+}) , \omega_{Y}^{\bullet}[-n] ) 
\xraq{\eqref{pseudo_can_map}^*} \RcHom_{\cO_Y}(  \cO(K_Y) , \omega_{Y}^{\bullet}[-n] ) 
\end{split}
\end{equation}
is equal to the composition
\begin{equation}  \label{compatible_lemma_1_map_3}
\cO_Y \xraq{\Delta} \RcHom_{\cO_Y}( \cO(K_Y) , \cO(K_Y) ) \xraq{\eqref{O_K_Y_omega_map}_*} \RcHom_{\cO_Y}( \cO(K_Y) , \omega_{Y}^{\bullet}[-n] )
\end{equation}
\elm

\bpf
Notice that a map $\cO_Y \ra \cF$ in $\cD(\QCoh(Y))$ is uniquely determined by its image of the global section $1 \in H^0(Y,\cO_Y)$.
It is clear that the composition \eqref{compatible_lemma_1_map_3} sends $1$ to \eqref{O_K_Y_omega_map}.
Hence, it suffices to show that the compositions \eqref{compatible_lemma_1_map_1} also sends $1$ to the same map.
By the description \eqref{local_adj_map_Phi_descr}, this statement is equivalent to the commutativity of the outermost square of
\begin{equation*}
\begin{tikzcd} [column sep = 40]
{\bm R}\pi^+_* \cO(K_{X^+}) \ar[r, "{\bm R}\pi^+_*(\Phi_{X^+})"]
& {\bm R}\pi^+_* \cH^{-n}(\omega_{X^+}^{\bullet}) \ar[r, "{\bm R}\pi^+_*(\iota)"]
& {\bm R}\pi^+_* \omega_{X^+}^{\bullet}[-n] \ar[dd, "\Tr_{X^+/Y}"]\\
\pi^+_* \cO(K_{X^+}) \ar[r, "\pi^+_*(\Phi_{X^+})"] \ar[u, "\iota"'] \ar[d, equal]
& \pi^+_*  \cH^{-n}(\omega_{X^+}^{\bullet}) \ar[d, "\!\,^0 \Tr_{X^+/Y}"] \ar[u, "\iota"'] & \\
\cO(K_Y)  \ar[r, "\Phi_{Y}"]
& \cH^{-n}(\omega_Y^{\bullet}) \ar[r, "\iota"]
& \omega_Y^{\bullet}[-n]
\end{tikzcd}
\end{equation*}
We verify the commutativity of each square, where the lower left one is given by \eqref{Phi_maps_commute}.
\epf

Notice that Lemma \ref{compatible_lemma_1} asserts that, although the composition \eqref{compatible_lemma_1_map_1} depends on $\pi^+ : X^+ \ra Y$, the resulting map can be described in terms of $\Phi_Y$, without any reference to $X^+$. Thus, the corresponding compositions for $X^-$ and $X^+$ coincide. This compatiblity between the local adjunction maps for $X^-$ and $X^+$ will be part of the condition for homological flips/flops to be introduced in the next section.

\brm  \label{CM_cond_rem}
We discuss some of the assumptions we have imposed. Starting from the discussion \eqref{pi_shriek_KX}, we have assumed that $X^-$ and $X^+$ are Cohen-Macaulay. 
Without this assumption, the maps \eqref{pi_shriek_KX}, and hence \eqref{local_adj_map_Phi}, still exists, but may not be an isomorphism. Since all the maps are pointing to the ``correct'' directions, Lemma \ref{compatible_lemma_1} also holds without the Cohen-Macaulay condition. 

Another condition we have imposed is \eqref{pseudo_can_cond}. Without this assumption, the map \eqref{pseudo_can_map} in the statement of Lemma \ref{compatible_lemma_1} cannot be defined.
\erm

%
%

%
%

\vspace{0.2cm}

Now we work in the following setting:
\begin{equation}  \label{crepant_situation}
\parbox{40em}{The diagram \eqref{Xm_Y_Xp_2} of birational contractions satisfies \eqref{pseudo_can_cond}, both $X^-$ and $X^+$ are Cohen-Macaulay, and both $\pi^-$ and $\pi^+$ are strongly crepant in the sense that $K_Y$ is Cartier, and $(\pi^{\pm})^* K_Y = K_{X^{\pm}}$.}
\end{equation}
This is satisfied for (1) strongly crepant flops between Cohen-Macaulay normal varieties; and (2) diagrams \eqref{Xm_Y_Xp_2} of birational contractions where both $X^-$ and $X^+$ are Calabi-Yau (\ie, both are Cohen-Macaulay and have trivial canonical divisor $K_{X^{\pm}} \equiv 0$).

In this case, there are canonical identifications%
\footnote{We write equality signs because they are equal as subsheaves of the sheaves of rational functions $\mathscr{K}_{X^-}$ and $\mathscr{K}_{X^+}$.}
\begin{equation*}  
\cO(K_{X^-}) \, = \, (\pi^-)^*\cO(K_Y) \qquad \text{and} \qquad 
\cO(K_{X^+}) \, = \, (\pi^+)^*\cO(K_Y)
\end{equation*}

Since we assume that $K_Y$ is Cartier, the object $\omega_Y^{\prime \bullet} := \omega_Y^{\bullet}[-n] \otimes_{\cO_Y} \cO(-K_Y)$ in $\Dbcoh(\QCoh(Y))$ is a dualizing complex.
Moreover, if we apply \eqref{shriek_tensor} to $\cF = \omega_Y^{\bullet}[-n]$ and $\cG = \cO(-K_Y)$,
we see that there is an isomorphism
\begin{equation}  \label{omega_K_to_pi_shriek}
\omega_{X^+}^{\bullet}[-n] \otimes_{\cO_{X^+}} \cO(-K_{X^+}) \xraq{\cong} (\pi^+)^!(\omega_Y^{\prime \bullet})
\end{equation}

On the other hand, $\Phi_{X^+}$ gives a canonical isomorphism
\begin{equation}  \label{flop_Phi_OX_pre} \cO_{X^+} = \cO(K_{X^+}) \otimes_{\cO_{X^+}} \cO(-K_{X^+})
\xraq[\cong]{ (\iota \circ \Phi_{X^+}) \otimes \id }  
\omega_{X^+}^{\bullet}[-n] \otimes_{\cO_{X^+}} \cO(-K_{X^+}) 
\end{equation}

Composing with these two maps, we obtain a canonical isomorphism
\begin{equation}  \label{flop_Phi_OX}
\widetilde{\Phi}'_{X^+} \, : \, \cO_{X^+} \xraq[\cong]{\eqref{flop_Phi_OX_pre}}
\omega_{X^+}^{\bullet}[-n] \otimes_{\cO_{X^+}} \cO(-K_{X^+})  \xraq[\cong]{\eqref{omega_K_to_pi_shriek}} (\pi^+)^!(\omega_Y^{\prime \bullet})
\end{equation}

We can again use this isomorphism to express Grothendieck duality for $\pi^+$ in terms of $\cO_{X^+}$ and $\omega_Y^{\prime \bullet}$. Namely, for any $\cF \in \Dmcoh(\QCoh(X^+))$, there is a local adjunction isomorphism
\begin{equation}  \label{local_adj_map_Phi_prime}
{\bm R}\pi^+_* \RcHom_{\cO_{X^+}}(\cF , \cO_{X^+}) \xraq{\cong} 
\RcHom_{\cO_Y}( {\bm R}\pi^+_* \cF , \omega_{Y}^{\prime \bullet} )
\end{equation}

Unravelling the definitions, we see the Lemma \ref{compatible_lemma_1} can be rewritten in this case in the following form, which will also be used in the next section in the proof of Theorem \ref{HFF_main_examples}:
\blm  \label{compatible_lemma_2}
The composition 
\begin{equation}  \label{compatible_lemma_2_map1}
\begin{split}
\cO_Y &\raq {\bm R}\pi^+_* \cO_{X^+} \, = \, 
{\bm R}\pi^+_* \RcHom_{\cO_{X^+}}(\cO_{X^+} , \cO_{X^+}) \\
&\xraq{ \eqref{local_adj_map_Phi_prime} } 
\RcHom_{\cO_Y}( {\bm R}\pi^+_* \cO_{X^+} , \omega_{Y}^{\prime \bullet} ) 
\raq \RcHom_{\cO_Y}(  \cO_Y , \omega_{Y}^{\prime \bullet} ) 
\, = \, \omega_{Y}^{\prime \bullet} 
\end{split}
\end{equation}
is equal to the composition
\begin{equation}  \label{compatible_lemma_2_map}
\cO_Y = \cO(K_Y) \otimes_{\cO_Y} \cO(-K_Y) \xraq{\eqref{O_K_Y_omega_map} \otimes \id}  \omega_{Y}^{\bullet}[-n]  \otimes_{\cO_Y} \cO(-K_Y) \, = \, \omega_{Y}^{\prime \bullet}
\end{equation}
\elm

\bpf
By definition, the map \eqref{omega_K_to_pi_shriek} is adjoint to the map
\begin{equation}  \label{omega_K_to_pi_shriek_adj}
\begin{split}
&{\bm R}\pi^+_*( \omega_{X^+}^{\bullet}[-n] \otimes_{\cO_{X^+}} \cO(-K_{X^+}) ) \, = \, {\bm R}\pi^+_*( \omega_{X^+}^{\bullet}[-n]) \otimes_{\cO_Y} \cO(-K_Y) \\ & \xraq{\Tr_{X^+/Y} \otimes \id} \omega_{Y}^{\bullet}[-n] \otimes_{\cO_Y} \cO(-K_Y) \, = \, \omega_Y^{\prime \bullet}
\end{split}
\end{equation}
In particular, for $\cF = \cO_{X^+}$, the map \eqref{local_adj_map_Phi_prime} sends the global section $1 \in \Hom_{\cD(\QCoh(X^+))}(\cO_{X^+},\cO_{X^+})$ to the global section
of $\RcHom_{\cO_Y}( {\bm R}\pi^+_* \cO_{X^+} , \omega_{Y}^{\prime \bullet} )$ given as the composition 
\begin{equation}  \label{Rpi_O_to_pi_shriek}
{\bm R}\pi^+_* \cO_{X^+} \xraq{{\bm R}\pi^+_*( \ref{flop_Phi_OX_pre} )} 
{\bm R}\pi^+_* (\omega_{X^+}^{\bullet}[-n] \otimes_{\cO_{X^+}} \cO(-K_{X^+}))  
\xraq{\eqref{omega_K_to_pi_shriek_adj}} \omega_Y^{\prime \bullet}
\end{equation}

Since a map from $\cO_Y$ to any other object in $\cD(\QCoh(Y))$ is uniquely determined by where it sends the global section $1$ to, this shows that \eqref{compatible_lemma_2_map1} can be rewritten as the left column of the following commutative diagram:
\begin{equation*}
\begin{tikzcd}
\cO_Y \ar[r, equal] \ar[d] & \cO_Y \ar[d] \\
{\bm R}\pi^+_* \cO_{X^+} \ar[r, equal] \ar[d, equal]
& {\bm R}\pi^+_* \cO_{X^+} \ar[d, "{\bm R}\pi^+_*(\Delta)"] \\
{\bm R}\pi^+_*( \cO(K_{X^+})  \otimes_{\cO_{X^+}} \cO(-K_{X^+}) ) \ar[r, "\cong"]
\ar[d, "{\bm R}\pi^+_*((\iota \circ \Phi_{X^+}) \otimes \id)"']
 &  {\bm R}\pi^+_* \RcHom_{\scO_{X^+}}( \cO(K_{X^+}) , \cO(K_{X^+}) ) \ar[d, "{\bm R}\pi^+_*(\iota \circ \Phi_{X^+})"] \\
{\bm R}\pi^+_*( \omega_{X^+}^{\bullet}[-n] \otimes_{\cO_{X^+}} \cO(-K_{X^+}) ) \ar[r, "\cong"] \ar[dd, "\eqref{omega_K_to_pi_shriek_adj}"'] &  {\bm R}\pi^+_* \RcHom_{\scO_{X^+}}( \cO(K_{X^+}) , \omega_{X^+}^{\bullet}[-n] ) \ar[d, "\cong"] \\
& \RcHom_{\scO_{Y}}( {\bm R}\pi^+_* \cO(K_{X^+}) , \omega_{Y}^{\bullet}[-n] ) \ar[d] \\
\omega_{Y}^{\bullet}[-n] \otimes_{\cO_{Y}} \cO(-K_{Y}) \ar[r, "\cong"]
& \RcHom_{\scO_{Y}}(  \cO(K_{Y}) , \omega_{Y}^{\bullet}[-n] )
\end{tikzcd}
\end{equation*}

Notice that the right column is the map \eqref{compatible_lemma_1_map_1}, so that these two maps coincide under the bottom isomorphism.
The result then follows from Lemma \ref{compatible_lemma_1} since the maps  \eqref{compatible_lemma_2_map} and \eqref{compatible_lemma_1_map_3} also coincide under the same bottom isomorphism. 
%
\epf

\brm
\begin{enumerate}
	\item The same remarks as in Remark \ref{CM_cond_rem} about the Cohen-Macaulay condition on $X^-$ and $X^+$ hold in the present context of \eqref{crepant_situation}.
	\item If $Y$ is Cohen-Macaulay (hence Gorenstein), then the dualizing complex $\omega_Y^{\prime \bullet}$ is isomorphic to $\cO_Y$ via the map $\Phi_Y$.
\end{enumerate}
\erm

\section{Homological flips and homological flops}  \label{HFF_sec}

In this section, we first extract some homological structures from some classes of flips and flops. These will motivate the definition of homological flip/flop in Definition \ref{HFF_def} below.
After that, we prove that these classes of flips and flops actually form homological flips/flops (see Theorem \ref{HFF_main_examples}). The crucial point is to verify a compatibility condition \eqref{HFF_compatible} across the flip/flop. Finally, we prove the main result (Theorem \ref{HFF_main_thm} and \ref{HFF_Gorenstein}) which concerns a duality of local cohomology groups.

Let $(Y,\cA)$ be a pair as in \eqref{sheaf_A_setting}. Fix a dualizing complex $\omega_Y^{\prime \prime \bullet} \in \Dbcoh(\QCoh(Y))$, and consider the involution \eqref{D_Y_def} defined using $\omega_Y^{\prime \prime \bullet}$.
By Theorem \ref{GGM_thm_nonaffine}, the object $(\pi^+)^! (\omega_Y^{\prime \prime \bullet})$ has a description in terms of the functor $\bD_Y$. Namely, there is a canonical isomorphism in $\cD(\Gr(\cA))$:
\begin{equation} \label{R_pi_shriek_D_Y}
\cR^+((\pi^+)^!( \omega_Y^{\prime \prime \bullet} )) \, \cong \, \bD_Y(\Ce_{\scI^+}(\cA))
\end{equation}

We will work with the following settings ({\it cf.} Definition \ref{flip_flop_def}):

\begin{equation}  \label{pair_flip_situation}
\parbox{40em}{Consider a flip \eqref{Xm_Y_Xp_2}, where both $X^-$ and $X^+$ are Cohen-Macaulay, and there exists $a > 0$ such that $K_{X^{\pm}} = a D^{\pm}$. Take $\cA := \bigoplus_{i \in \bZ} \cO_Y(iD_Y)$.}
\end{equation}
\begin{equation} \label{pair_flop_situation}
\parbox{40em}{Consider a log flip \eqref{Xm_Y_Xp_2}, where both $X^-$ and $X^+$ are Cohen-Macaulay, and $K_Y$ is Cartier. Take $\cA := \bigoplus_{i \in \bZ} \cO_Y(iD_Y)$.}
\end{equation}

In the situation \eqref{pair_flip_situation}, there are isomorphisms
\begin{equation}  \label{A_tilde_to_omega_flip}
\widetilde{\cA(a)}_{X^{\pm}}
\xraq[\cong]{\eqref{O_iD_to_A_tilde}} 
\cO_{X^{\pm}}(K_{X^{\pm}}) 
\xraq[\cong]{\eqref{pi_shriek_KX}} (\pi^{\pm})^!(\omega_Y^{\bullet}[-n])
\qquad \text{ in }  \Dbcoh(\QCoh(X^{\pm}))
\end{equation} 

In the situation \eqref{pair_flop_situation}, there are isomorphisms
\begin{equation} \label{A_tilde_to_omega_flop}
\widetilde{\cA}_{X^{\pm}}
\, = \, 
\cO_{X^{\pm}}
\xraq[\cong]{\eqref{flop_Phi_OX}}
\pi^!(\omega_Y^{\prime \bullet})
\qquad \text{ in } \Dbcoh(\QCoh(X^{\pm}))
\end{equation} 

Both of the equations \eqref{A_tilde_to_omega_flip} and \eqref{A_tilde_to_omega_flop} can be written in the same way
\begin{equation}  \label{A_a_tilde_to_pi_shriek}
\widetilde{\cA(a)}_{X^{\pm}} \xraq{\cong} (\pi^{\pm})^!(\omega_Y^{\prime \prime \bullet})
\end{equation}
where we take $\omega_Y^{\prime \prime \bullet} := \omega_Y^{\bullet}[-n]$ and $a>0$ in the case \eqref{pair_flip_situation}; and $\omega_Y^{\prime \prime \bullet} := \omega_Y^{\prime \bullet}$ and $a=0$ in the case \eqref{pair_flop_situation}.
Notice that $\omega_Y^{\prime \prime \bullet}$ is a dualizing complex in both cases.

For notational simplicity, we focus on the $X^+$-side in the following discussion. Applying  \eqref{Cech_R_functor_nonaffine} to $\cM = \cA(a)$, we see that there is a canonical map 
\begin{equation}  \label{Cech_R_pi_shriek}
\Ce_{\scI^+}(\cA)(a) \ra \cR^+( \widetilde{\cA(a)}_{X^+}) \xra{\cong} \cR^+( (\pi^+)^!( \omega_Y^{\prime \prime \bullet}))
\end{equation}
where the weight $i$ component of the first map
is given by 
\begin{equation}  \label{Cech_R_pi_shriek_comp}
{\bm R}\pi^+_* ( \widetilde{\cA(a+i)})  \raq {\bm R}\pi^+_* \RcHom_{\cO_{X^+}}(\widetilde{\cA(-i)}, \widetilde{\cA(a)})
\end{equation} by the description 
\eqref{Cech_R_functor_nonaffine_comp}. Thus, if $\cA$ is $\tfrac{1}{d}$-Cartier, then  the map \eqref{Cech_R_pi_shriek} is a quasi-isomorphism in weight components $i \in d \bZ$. 
In fact, it is often a quasi-isomorphism for all weight components. 
Namely, if the pair $(Y,\cA)$ determines a Cohen-Macaulay log flip as in Definition \ref{determine_log_flip}, then the map \eqref{Cech_R_pi_shriek_comp} can be rewritten as the derived pushforward ${\bm R}\pi^+_*(-)$ of the map \eqref{O_D_RHom} below, for $X := X^+$, $D := D^+$ and $D' := aD^+$, and hence is a quasi-isomorphism.


\blm
Suppose that $X$ is a  Cohen-Macaulay normal variety over $k$, and $D \in \WDiv(X)$ is a Weil divisor on $X$ such that the sheaf $\cO_X(-iD)$ is maximal Cohen-Macaulay. If $D' \in \WDiv(X)$ is such that $\cO_X(D')$ is a dualizing complex, then the canonical map 
\begin{equation}  \label{O_D_RHom}
\cO_X(D' + iD) \raq \RcHom_{\cO_X}(\cO_X(-iD), \cO_X(D'))
\end{equation}
is a quasi-isomorphism.
\elm

\bpf
The $0$-th cohomology sheaves of \eqref{O_D_RHom} is the map $\cO_X(D' + iD) \ra \cHom_{\cO_X}(\cO_X(-iD), \cO_X(D'))$, which is always an isomorphism. Thus, it suffices to show that the complex $\RcHom_{\cO_X}(\cO_X(-iD), \cO_X(D'))$ has no higher cohomology sheaves. Since $\cO_X(D')$ is assumed to be a dualizing complex, this is equivalent to the sheaf $\cO_X(-iD)$ being maximal Cohen-Macaulay (see, e.g., \cite[Tag 0A7U]{Sta}).
\epf

Combining \eqref{R_pi_shriek_D_Y} and \eqref{Cech_R_pi_shriek}, we obtain a map
\begin{equation} \label{Phi_plus}
\Phi^+ \, : \, \Ce_{\scI^+}(\cA)(a) \raq \bD_Y(\Ce_{\scI^+}(\cA)) 
\end{equation}
in $\cD(\Gr(\cA))$, which is an isomorphism in weight $i$ if $\cO_{X^+}(-i D^+)$ is maximally Cohen-Macaulay.

Simiarly discussion holds for $X^-$ and $\scI^-$ in place of $X^+$ and $\scI^+$, and we likewise have a map 
\begin{equation} \label{Phi_minus}
\Phi^- \, : \, \Ce_{\scI^-}(\cA)(a) \raq \bD_Y(\Ce_{\scI^-}(\cA))  
\end{equation}
in $\cD(\Gr(\cA))$, which is an isomorphism in weight $i$ if $\cO_{X^-}(-i D^-)$ is maximally Cohen-Macaulay.
 Moreover, as we will see in Theorem \ref{HFF_main_examples} below, the maps \eqref{Phi_plus} and \eqref{Phi_minus} thus obtained satisfy a certain compatibility condition, as formalized in the following
\bdf  \label{HFF_def}
A \emph{weak homological flip} (resp. \emph{weak homological flop}) consists of a sextuple $(Y,\omega_Y^{\prime \prime \bullet}, \cA, a, \Phi^-,\Phi^+)$ where
\begin{enumerate}
	\item $Y$ is a Noetherian separated scheme with a dualizing complex $\omega_Y^{\prime \prime \bullet} \in \Dbcoh(\QCoh(Y))$.
	\item $\cA$ is a quasi-coherent sheaf of Noetherian $\bZ$-graded rings on $\cA$ such that $\cA_0$ is coherent over $Y$.
	\item $a > 0$ (resp. $a=0$) is an integer
	\item $\Phi^-$ and $\Phi^+$ are maps in $\cD(\Gr(\cA))$ 
	\begin{equation*}  
	\begin{split}
	\Phi^+ \, &: \, \Ce_{\scI^+}(\cA)(a) \raq \bD_Y(\Ce_{\scI^+}(\cA)) \\
	\Phi^- \, &: \, \Ce_{\scI^-}(\cA)(a) \raq \bD_Y(\Ce_{\scI^-}(\cA))
	\end{split}
	\end{equation*}
	that are quasi-isomorphisms in weights $i \in d\bZ$ if $\cA$ is $\tfrac{1}{d}$-Cartier. Here $\bD_Y$ is the functor \eqref{D_Y_def} defined using $\omega_Y^{\prime \prime \bullet}$.
\end{enumerate}
The maps $\Phi^-$ and $\Phi^+$ are required to be compatible in the sense that the diagram 
\begin{equation}  \label{HFF_compatible}
\begin{tikzcd} [row sep = 0]
& \Ce_{\scI^+}(\cA)(a) \ar[r, "\Phi^+"] & \bD_Y( \Ce_{\scI^+}(\cA) ) \ar[rd, "\bD_Y(\eta^+)"] & \\
\cA(a) \ar[ru, "\eta^+"] \ar[rd, "\eta^-"']  & & & \bD_Y(\cA) \\
& \Ce_{\scI^-}(\cA)(a) \ar[r, "\Phi^-"] & \bD_Y( \Ce_{\scI^-}(\cA) ) \ar[ru, "\bD_Y(\eta^-)"'] &
\end{tikzcd}
\end{equation}
commutes in $\cD(\Gr(\cA))$, where $\eta^{\pm} : \cA \ra \Ce_{\scI^{\pm}}(\cA)$ are the adjunction units as in \eqref{RGam_Ce_seq}. 

A weak homological flip/flop is said to be a \emph{homological flip/flop} if the maps $\Phi^-$ and $\Phi^+$ are quasi-isomorphisms in all weight degrees.
\edf

We also introduce the following extra conditions:

\bdf  \label{canonical_range_def}
A (weak) homological flip/flop $(Y,\omega_Y^{\prime \prime \bullet}, \cA, a, \Phi^-,\Phi^+)$ is said to have \emph{positive canonical vanishing} if $\RGam_{\scI^+}(\cA)_i \simeq 0$ for all $i \geq a$. It is said to have \emph{negative canonical vanishing} if $\RGam_{\scI^-}(\cA)_i \simeq 0$ for all $i \leq a$.
It is said to be \emph{canonical vanishing} if it satisfies both.
\edf




Our discussion above leads to the following first main result of this section, which gives examples of (weak) homological flip/flop:
\bthm  \label{HFF_main_examples}
\begin{enumerate}
	\item In the situation \eqref{pair_flip_situation}, take $\omega_Y^{\prime \prime \bullet} := \omega_Y^{\bullet}[-n]$, and let $\Phi^-$ and $\Phi^+$ be the maps \eqref{Phi_minus} and \eqref{Phi_plus}, then the sextuple $(Y,\omega_Y^{\prime \prime \bullet}, \cA, a, \Phi^-,\Phi^+)$ is a weak homological flip. If each $\cO_{X^{\pm}}(iD^{\pm})$ is maximal Cohen-Macaulay, then $(Y,\omega_Y^{\prime \prime \bullet}, \cA, a, \Phi^-,\Phi^+)$ is a homological flip.
	
	\item In the situation \eqref{pair_flop_situation}, take $\omega_Y^{\prime \prime \bullet} := \omega_Y^{\prime \bullet}$ (see the discussion following \eqref{crepant_situation}), and let $\Phi^-$ and $\Phi^+$ be the maps \eqref{Phi_minus} and \eqref{Phi_plus} for $a = 0$, then the sextuple $(Y,\omega_Y^{\prime \prime \bullet}, \cA, 0, \Phi^-,\Phi^+)$ is a weak homological flop. If each $\cO_{X^{\pm}}(iD^{\pm})$ is maximal Cohen-Macaulay, then $(Y,\omega_Y^{\prime \prime \bullet}, \cA, 0, \Phi^-,\Phi^+)$ is a homological flop.
\end{enumerate}


In either case, assume furthermore that $X^-$ and $X^+$ are projective over the base field (assumed to be of characteristic zero). Then
\begin{enumerate}
	\item[(a)] If $X^+$ has rational singularities and $D^+$ is Cartier, then the (weak) homological flip/flop has positive canonical vanishing. If $Y$ also has rational singularities, then $\RGam_{\scI^+}(\cA)_0 \simeq 0$ as well.
	\item[(b)] If $X^-$ has rational singularities and $D^-$ is Cartier, then the (weak) homological flip/flop has negative canonical vanishing, and $Y$ has rational singularities.
\end{enumerate}
\ethm

\bpf
In the discussion preceeding Definition \ref{HFF_def}, we have already proved all the statements in both cases (1) and (2), except the commutativity of \eqref{HFF_compatible} in $\cD(\Gr(\cA))$. For any $\cM \in \cD(\Gr(\cA))$, we have 
$\Hom_{\cD(\Gr(\cA))}(\cA(a),\cM) \cong \bH^0(Y, \cM_{-a})$, so that it suffices to show that the global section $1 \in \bH^0(Y,\cA(a)_{-a})$ is sent to the same element in $\bH^0( Y, \bD_Y(\cA)_{-a}) = \Hom_{\cD(\QCoh(Y))}( \cA_a , \omega_Y^{\prime \prime \bullet} )$ under the two routes in \eqref{HFF_compatible}. 

Consider the upper route of \eqref{HFF_compatible}, whose definition involves the maps $\eta^+ : \cA \ra \Ce_{\scI^+}(\cA)$, as well as the maps \eqref{Cech_R_functor_nonaffine} and \eqref{cR_pi_shriek_isom_sheaf}. 
In Section \ref{DGrA_sec}, we have carefully kept track of each weight degrees of these maps, in the form of the commutative diagrams \eqref{eta_M_weight_description}, \eqref{Cech_R_functor_nonaffine_comp} and \eqref{cR_pi_shriek_isom_2}. 
Thus, if we unravel the definitions, we see that the upper route of \eqref{HFF_compatible} at weight $-a$ is given by
\begin{equation}  \label{HFF_comp_weight_minus_a}
\begin{split}
& \cO_Y \raq {\bm R}\pi^+_* \cO_{X^+} \xraq{ {\bm R}\pi^+_*(\Delta) }
{\bm R}\pi^+_* \RcHom_{\cO_{X^+}}(\widetilde{\cA(a)}_{X^+} , \widetilde{\cA(a)}_{X^+}) \\
& \xraq{\eqref{A_a_tilde_to_pi_shriek}}  {\bm R}\pi^+_* \RcHom_{\cO_{X^+}}( \widetilde{\cA(a)}_{X^+} ,  (\pi^+)^!(\omega_Y^{\prime \prime \bullet}) ) 
\, \cong \, 
\RcHom_{\cO_Y}( {\bm R}\pi^+_* (\widetilde{\cA(a)}_{X^+}) , \omega_Y^{\prime \prime \bullet}) ) \\
& \xraq{\iota^*} \RcHom_{\cO_Y}(  \pi^+_* (\widetilde{\cA(a)}_{X^+}) , \omega_Y^{\prime \prime \bullet}) ) 
\xraq{\eqref{M_0_pi_M}^*} \RcHom_{\cO_Y}(  \cA_a , \omega_Y^{\prime \prime \bullet}) ) 
\end{split}
\end{equation}

For the case (2) of flops, with $a = 0$, we have $\widetilde{\cA(a)}_{X^+} = \cO_{X^+}$, so that the composition \eqref{HFF_comp_weight_minus_a} is precisely \eqref{compatible_lemma_2_map1}, and hence is equal to \eqref{compatible_lemma_2_map} by Lemma \ref{compatible_lemma_2}. Since the map \eqref{compatible_lemma_2_map} is defined purely in terms of structures on $Y$, it is the same for $X^+$ and $X^-$.

For the case (1) of flips, the map \eqref{A_a_tilde_to_pi_shriek} is given by \eqref{A_tilde_to_omega_flip}, which factors through the isomorphism \eqref{O_iD_to_A_tilde}. Thus, the part ``$\eqref{A_a_tilde_to_pi_shriek}\circ {\bm R}\pi^+_*(\Delta) $" appearing in \eqref{HFF_comp_weight_minus_a} may be rewritten as
\begin{equation*}
{\bm R}\pi^+_* \cO_{X^+} \xraq{ \eqref{psi_A_tilde_O_iD} }
{\bm R}\pi^+_* \RcHom_{\cO_{X^+}}(\widetilde{\cA(a)}_{X^+} , \cO(K_{X^+}))  \xraq{\eqref{pi_shriek_KX}}  {\bm R}\pi^+_* \RcHom_{\cO_{X^+}}( \widetilde{\cA(a)}_{X^+} ,  (\pi^+)^!(\omega_Y^{\bullet}[-n]) ) 
\end{equation*}

Replace it as such, \eqref{HFF_comp_weight_minus_a} then becomes the composition of $\cO_Y \ra {\bm R} \pi^+_* \cO_{X^+}$ with the right column of the commutative diagram
\begin{equation*}
\begin{tikzcd}
{\bm R}\pi^+_* \cO_{X^+} \ar[r, equal] \ar[d, "{\bm R}\pi^+_*(\Delta)"'] & {\bm R}\pi^+_* \cO_{X^+} \ar[d, "\eqref{psi_A_tilde_O_iD}"]\\
{\bm R}\pi^+_* \RcHom_{\cO_{X^+}}(\cO(K_{X^+}) , \cO(K_{X^+})) \ar[r, "\eqref{psi_A_tilde_O_iD}^*"] \ar[d, "\eqref{pi_shriek_KX}"']
& {\bm R}\pi^+_* \RcHom_{\cO_{X^+}}(\widetilde{\cA(a)}_{X^+} , \cO(K_{X^+})) \ar[d, "\eqref{pi_shriek_KX}"] \\
{\bm R}\pi^+_* \RcHom_{\cO_{X^+}}(  \cO(K_{X^+}) ,  (\pi^+)^!(\omega_Y^{\bullet}[-n]) ) \ar[r, "\eqref{psi_A_tilde_O_iD}^*"]
\ar[d, "\cong"'] & {\bm R}\pi^+_* \RcHom_{\cO_{X^+}}( \widetilde{\cA(a)}_{X^+} ,  (\pi^+)^!(\omega_Y^{\bullet}[-n]) ) \ar[d, "\cong"] \\
\RcHom_{\cO_Y}( {\bm R}\pi^+_* (\cO(K_{X^+}) , \omega_Y^{\bullet}[-n]) ) \ar[d, "\iota^*"'] \ar[r, "\eqref{psi_A_tilde_O_iD}^*"]
& \RcHom_{\cO_Y}( {\bm R}\pi^+_* (\widetilde{\cA(a)}_{X^+}) , \omega_Y^{\bullet}[-n]) ) \ar[d, "\iota^*"]\\
\RcHom_{\cO_Y}(  \pi^+_* (\cO(K_{X^+})) , \omega_Y^{\bullet}[-n]) ) \ar[d, equal] \ar[r, "\eqref{psi_A_tilde_O_iD}^*"]
& \RcHom_{\cO_Y}(  \pi^+_* (\widetilde{\cA(a)}_{X^+}) , \omega_Y^{\bullet}[-n]) ) \ar[d, "\eqref{M_0_pi_M}^*"] \\
\RcHom_{\cO_Y}(  \cO(K_Y) , \omega_Y^{ \bullet}[-n]) ) \ar[r, equal]
& \RcHom_{\cO_Y}(  \cA_a , \omega_Y^{ \bullet}[-n]) ) 
\end{tikzcd}
\end{equation*}
where the bottom square commutes because of \eqref{O_iD_to_A_tilde_comp_prop}. Moreover, the left column, when precomposed with $\cO_Y \ra {\bm R} \pi^+_* \cO_{X^+}$, is precisely the map \eqref{compatible_lemma_1_map_1}. By Lemma \ref{compatible_lemma_1}, it is therefore the same for $X^+$ and $X^-$. This shows that \eqref{HFF_compatible} commutes in this case.

For the statements (a) and (b), notice that the long exact sequence associated to \eqref{RGam_Ce_seq} for $\cM = \cA$ gives at weight $i$ the exact sequence
\begin{equation*}
0 \raq \cH^0(\RGam_{\scI^{\pm}}(\cA)_i) \xraq{\epsilon} \cA_i \xraq[\cong]{\eqref{M_0_pi_M}} \pi^{\pm}_*(\widetilde{\cA(i)}_{X^{\pm}}) \xraq{\delta} \cH^1(\RGam_{\scI^{\pm}}(\cA)_i) \raq 0
\end{equation*}
and the isomorphisms $ {\bm R}^j \pi^{\pm}_*(\widetilde{\cA(i)}_{X^{\pm}}) \xraq{\delta} \cH^{j+1}(\RGam_{\scI^{\pm}}(\cA)_i)$ for any $j \geq 1$.
This shows that $\cH^{\leq 1}\RGam_{\scI^{\pm}}(A)_i = 0$, and we have
\begin{equation*}
\RGam_{\scI^{\pm}}(A)_i \simeq 0 \quad \text{ if and only if } \quad {\bm R}^{j}\pi^{\pm}_*( \widetilde{\cA(i)}_{X^{\pm}} ) \text{ for all } j > 0
\end{equation*}

Notice also that in both cases, the sheaf $\widetilde{\cA(a)}_{X^{\pm}}$ is linearly $\pi^{\pm}$-equivalent to the canonical sheaf. Namely, in the case (1), we have $\widetilde{\cA(a)}_{X^{\pm}} = \cO(K_{X^{\pm}})$; while in the case (2), we have 
$\widetilde{\cA}_{X^{\pm}} = \cO_{X^{\pm}} \cong \cO(K_{X^{\pm}}) \otimes (\pi^{\pm})^*\cO(-K_Y)$.
Thus, both statements about positive/negative canonical vanishing follow from the Kawamata-Viehweg vanishing Theorem in the form of Corollary \ref{KV_vanishing_2}. Moreover, if both $X^+$ and $Y$ have rational singularities, then ${\bm R}^j\pi^+_* \cO_{X^+} = 0$ for all $j > 0$, and hence $\RGam_{\scI^+}(\cA)_0 \simeq 0$. On the other hand, if $X^-$ has rational singularities, then Corollary \ref{KV_vanishing_2} implies in particular that
${\bm R}^j\pi^-_* \cO_{X^-} = 0$ for all $j > 0$, so that $Y$ has rational singularities.
%
\epf

\brm  \label{KV_klt_remark}
In the part about canonical vanishing, we have imposed the condition that $D^+$ (resp. $D^-$) is Cartier. This condition can be removed if we use a stronger form of the Kawamata-Viehweg vanishing Theorem (see, \eg, \cite[Theorem 2.70]{KM98} or \cite[Theorem 2.17.3]{Kol97}). For this, it suffices to assume that $X^{\pm}$ are proper with log terminal singularities (\ie, $(X^{\pm},0)$ are klt). Notice that in this case, \cite[Corollary 5.25]{KM98} shows that $\cO_{X^{\pm}}(iD^{\pm})$ is maximal Cohen-Macaulay for each $i \in \bZ$, a condition we required in parts (1), (2) of Theorem \ref{HFF_main_examples}. 
\erm

The following is the second main result of this section:

\bthm  \label{HFF_main_thm}
Let $(Y,\omega_Y^{\prime \prime \bullet}, \cA, a, \Phi^-,\Phi^+)$ be a 
homological flip/flop. Then there is a map $\Psi : \RGam_{I^+}(A)(a)[1] \ra \bD_Y( \RGam_{I^-}(A) ) $ in $\cD(\GrA)$ such that, if $c^-, c^+ \in \bZ$ are integers as in Lemma \ref{local_cohom_weight_bounded} for $\cM=\cA$, then $\Psi$ is a quasi-isomorphism in weights
\begin{equation} \label{HFF_main_thm_i_range}
 i \geq \max\{ -c^- , c^+ - a\}  \quad \text{or} \quad 
i \leq \min\{ -c^+,c^- -a \} 
\end{equation}
%
\ethm

\bpf
Consider the diagram
\begin{equation}  \label{HFF_big_diag_1}
\begin{tikzcd}
\cA(a) \ar[r, "\eta^+"] \ar[d, "\eta^-"] & \Ce_{\scI^+}(\cA)(a) \ar[d, "\bD_Y(\eta^+) \circ \Phi^+"] \ar[r, "\delta^+"] 
& \RGam_{\scI^+}(\cA)(a)[1] \ar[d, dashed, "\Psi"] \\
\Ce_{\scI^-}(\cA)(a) \ar[r, "\bD_Y(\eta^-) \circ \Phi^-"]  \ar[d, "\delta^-"]
& \bD_Y(\cA) \ar[r, "\bD_Y(\epsilon^-)"] \ar[d, "\bD_Y(\epsilon^+)"] 
& \bD_Y(\RGam_{\scI^-}(\cA)) \ar[d, dashed]\\
\RGam_{\scI^-}(\cA)(a)[1] \ar[r, dashed, "\Psi'"] &  \bD_Y(\RGam_{\scI^+}(\cA)) \ar[r, dashed] & Z
\end{tikzcd}
\end{equation}

Since $\Phi^-$ and $\Phi^+$ are isomorphisms, the first two rows and the first two columns are parts of distinguished triangles (take \eqref{RGam_Ce_seq} and their dual $\bD_Y(-)$). 
Moreover, the top left square is commutative by the compatbility condition \eqref{HFF_compatible} of a homological flip/flop.
Thus, by the $3\times 3$-lemma of a triangulated category (see, e.g., \cite[Proposition 1.1.11]{BBD82} or \cite[Lemma 2.6]{May01}), the object $Z$ and the dotted arrows $\Psi, \Psi'$ in diagram \eqref{HFF_big_diag_1} exist, making each row and column part of a distinguished triangle.

The cohomology of $\RGam_{\scI^+}(\cA)$ is zero in weight $\geq c^+$, while the cohomology of $\RGam_{\scI^-}(\cA)$ is zero in weight $\leq c^-$. Thus, viewing $Z$ as the cone of $\Psi$, we see that the cohomology of $Z$ is zero in weight $i \geq \max \{  -c^- , c^+ - a \}$.
Similarly, viewing $Z$ as the cone of $\Psi'$, we see that the cohomology of $Z$ is zero in weight $i \leq \min \{  -c^+ , c^- - a \}$. This completes the proof.
%
%
\epf

\brm
If $(Y,\omega_Y^{\prime \prime \bullet}, \cA, a, \Phi^-,\Phi^+)$ is a $\tfrac{1}{d}$-Cartier weak
homological flip/flop, then by a similar proof, one also obtains a map $\Psi : \RGam_{I^+}(A)(a)[1] \ra \bD_Y( \RGam_{I^-}(A) ) $ that is a quasi-isomorphism in weight $i \in d\bZ$ satisfying \eqref{HFF_main_thm_i_range}.
\erm

As a consequence, we give the following

%
%
%

\bpf[Proof of Theorem \ref{HFF_Gorenstein}]
Both statements concerning $\Psi$ follows directly from Theorem \ref{HFF_main_thm}, so that it suffices to prove that $A = \cA(U)$ is Gorenstein when $\Psi$ is an isomorphism.
%
Recall that $\bD_Y$ is local on the base $Y$ when applied to the object $\RGam_{\scI^{\pm}}(\cA) \in \cD_{\lc}^b(\GrA)$.
Thus, the isomorphism $\Psi$ restricts to one on any open subscheme $U \subset Y$, so we may assume that $Y=\Spec \, R$ is affine.

By the exact triangle \eqref{RGam_Ce_seq}, it suffices to show that $\RGam_{I^+}(A)$ and $\Ce_{I^+}(A)$ have finite injective dimension%
\footnote{It follows from \cite[Proposition 3.6.6]{BH93} that, for any $M \in \GrA$, we have ${\rm inj \, dim}_{\GrA}(M) \, \leq \, {\rm inj \, dim}_{\Mod(A)}(M) \, \leq \, {\rm inj \, dim}_{\GrA}(M) + 1$. Hence, it is the same whether finiteness of injective dimension is considered over $\GrA$ or $\Mod(A)$.}. By definition of a homological flip/flop, there is an isomorphism $\Phi^+ : \Ce_{I^+}(A)(a) \xra{\cong} \bD_Y(\Ce_{I^+}(A))$ in $\cD(\GrA)$. 
Since both $\Ce_{I^+}(A)$ and $\RGam_{I^-}(A)$ are represented by a finite complex of flat graded modules, the sought for statement follows from the following simple fact, applied to $N = \bD_Y(A)$, and $M$ being $\RGam_{I^-}(A)$ or $\Ce_{I^+}(A)$:
\begin{equation*}
\parbox{44em}{Suppose that $M, N \in \cD(\GrA)$ are objects such that $M$ has finite Tor-dimension, and $N$ has finite injective dimension, then $\RHomcom_A(M,N) \in \cD(\GrA)$ has finite injective dimension.}
\end{equation*}
Indeed, this follows from the adjunction $\RHomcom_A(L,\RHomcom_A(M,N)) \cong \RHomcom_A( L \otimes_A^{{\bm L}} M , N )$.
\epf

\brm  \label{flop_conn_min_rem}
If $X$ and $X'$ are two birational Gorenstein Calabi-Yau normal projective varieties with at most $\bQ$-factorial terminal singularities, then a result \cite{Kaw08} of Kawamata shows that they are connected by a finite sequence of flops.
All the intermediary varieties are still Gorenstein Calabi-Yau with at most $\bQ$-factorial terminal singularities, and the contracted variety $Y$ in each flop $X^- \ra Y \la X^+$ that appears is Gorenstein. Thus, by Theorem \ref{HFF_main_examples}, $X$ and $X'$ are connected by finitely many homological flops with canonical vanishing. Specifically, they satisfy the assumptions of Theorem \ref{main_thm_intro}(2).
\erm

\section{Derived category under flips and flops}  \label{Dcat_flip_flop_sec}

The major application of our theory of homological flip/flop is to relate the derived categories under flips and flops.
We will focus on the case of affine base $Y = \Spec \, R$. All results have formal extensions to the non-affine case ({\it cf.} Corollary \ref{Dperf_Ce_ff_flop} and \ref{Dperf_Ce_ff_flop_intro}).
By Theorem \ref{HFF_main_examples} and \ref{HFF_Gorenstein}, we see that a large class of flips and flops are controlled by (sheaves of) Noetherian $\bZ$-graded rings $A$ satisfying
the following condition, where  we set $a=0$ for flops, and $a=1$ for flips:
\begin{equation}  \label{Gor_and_local_cohom_dual_2}
\parbox{40em}{(i) $A$ is Gorenstein. \\
	(ii) There is an isomorphism $\Psi : \RGam_{I^+}(A)(a)[1] \xra{\cong }\bD_Y(\RGam_{I^-}(A))$ in $\cD(\GrA)$.}
\end{equation}

As was alluded to in the introduction (see the discussion following Theorem \ref{main_thm_intro}), this duality between $\RGam_{I^+}(A)$ and $\RGam_{I^-}(A)$ should allow one to relate the derived categories $\cD_{\IpTR}(\GrA)$ and $\cD_{\ImTR}(\GrA)$, which are closely related to the derived categories of $X^+$ and $X^-$ respectively (see \eqref{D_QCoh_ITR_intro}).
In this section, we give a preliminary instance of this relation.
A more sophisticated version is given in the paper \cite{Yeu20b}.

For any $M , N \in \cD(\GrA)$, notice that the fiber (\ie, cocone) of the canonical map
\begin{equation*}
\RHomcom_A(M , N ) \raq  \RHomcom_A( \Ce_{I^+}(M) , \Ce_{I^+}(N))
\end{equation*}
is given by $\RHomcom_A(M ,\RGam_{I^+}(N))$. As a result, we have the following
\blm  \label{CI_ff_RHom_RGam}
The functor $\Ce_{I^+} : \cD(\GrA) \ra \cD_{\IpTR}(\GrA)$ is fully faithful on the pair $(M , N[i])$ for all $i \in \bZ$ if and only if 
$\RHom_A(M ,\RGam_{I^+}(N)) := \RHomcom_A(M ,\RGam_{I^+}(N))_0 \simeq 0$.
\elm

%
\bdf  \label{partial_sec_def}
A full triangulated subcategory $\cE \subset \cD(\GrA)$ is said to be a \emph{positive} (resp. \emph{negative}) \emph{slice} if $\RHomcom_A(M ,\RGam_{I^+}(N))_0 \simeq 0$ (resp. $\RHomcom_A(M ,\RGam_{I^-}(N))_0 \simeq 0$) for all $M,N \in \cE$.
By Lemma \ref{CI_ff_RHom_RGam}, this holds if and only if the functor $\Ce_{I^{\pm}} : \cE \ra \cD_{I^{\pm}\text{-triv}}(\GrA)$ is fully faithful.

In this case, let $\cP^{\pm} \subset \cD_{I^{\pm}\text{-triv}}(\GrA)$ be the essential image of $\cE$ under $\Ce_{I^{\pm}}$, then $\cE$ is said to be a \emph{positive} (resp. \emph{negative}) \emph{slice} of $\cP^{\pm}$. Thus, $\cE$ is equivalent to $\cP^{\pm}$ in this case.
\edf

The following two results show that in the case of homological flops with canonical vanishing, if we restrict to $\Dperf(\GrA)$, then being a positive slice is equivalent to being a negative slice.

\bpp  \label{RHom_RGam_swap}
Assume that the Noetherian $\bZ$-graded ring $A$ satisfies \eqref{Gor_and_local_cohom_dual_2}(ii) for $a = 0$. Suppose that $M \in \Dperf(\GrA)$ and $N \in \cD(\GrA)$ is such that $\RGam_{I^+}(N) \in \cD_{\lc}(\GrA)$ (this holds, for example, if $N \in \cD_{\coh}(\GrA)$), then we have 
\begin{equation*}
\RHom_A(M,\RGam_{I^+}(N)) \simeq 0 \quad \Leftrightarrow \quad 
\RHom_A(N,\RGam_{I^-}(M)) \simeq 0
\end{equation*}
\epp

\bpf
We have the following sequence of isomorphisms in $\cD(\GrA)$:
\begin{equation*}
\begin{split}
\RHomcom_A(N,\RGam_{I^-}(M)) 
\, & = \,  \RHomcom_A(N,M \otimes_A^{{\bm L}} \RGam_{I^-}(A)) 
\\
&\cong \, \RHomcom_A(N \otimes_A^{{\bm L}} M^{\vee}, \RGam_{I^-}(A))  \\
&\cong \, \bD_Y( N \otimes_A^{{\bm L}} M^{\vee} \otimes_A^{{\bm L}} \RGam_{I^+}(A) ) [-1]\\
&\cong \, \bD_Y( \RHomcom_A(M, \RGam_{I^+}(N))) [-1]
\end{split}
\end{equation*}
Since $\bD_Y$ is involutive on $\cD_{\lc}(\GrA)$, the statement follows by taking the weight $0$ part of the above isomorphism. 
\epf

\brm
The isomorphism $\RHomcom_A(N,\RGam_{I^-}(M)) \cong \bD_Y( \RHomcom_A(M, \RGam_{I^+}(N))) [-1]$ in the proof of Proposition \ref{RHom_RGam_swap} also holds under a different set of assumptions on $M$ and $N$. For example, if we assume both (i) and (ii) of \eqref{Gor_and_local_cohom_dual_2}, then it holds for $M \in \Dbcoh(\GrA)$ and $N \in \cD(\GrA)$ such that $\RGam_{I^-}(N)$ has finite Tor-dimension.  
\erm

\bcor  \label{Dperf_Ce_ff_flop}
Assume that the Noetherian $\bZ$-graded ring $A$ satisfies \eqref{Gor_and_local_cohom_dual_2}(ii) for $a = 0$. Then for any full triangulated subcategory $\cE \subset \Dperf(\GrA)$, the functor $\Ce_{I^+} : \cE \ra \cD_{\IpTR}(\GrA)$ is fully faithful if and only if the functor $\Ce_{I^-} : \cE \ra \cD_{\ImTR}(\GrA)$ is fully faithful.
\ecor

\bcor  \label{Dperf_Ce_ff_flop_v2}
Denote by $\cD_c \subset \cD$ the subcategory of compact objects.
Suppose that both of the localization
functors $\Ce_{I^{\pm}} : \Dperf(\GrA) \ra \cD_{I^{\pm}\text{-triv}}(\GrA)_c$ admit fully faithful sections, then there are fully faithful exact functors in both directions between $\cD_{\IpTR}(\GrA)_c$ and $\cD_{\ImTR}(\GrA)_c$
\ecor


Let $\cD_{<w}(\GrA)$ be the full subcategory of $\cD(\GrA)$ consisting of objects $M \in \cD(\GrA)$ such that $M_i \simeq 0$ for all $i \geq w$. 
Notice that we often have $\RHomcom_A(M ,\RGam_{I^+}(N)) \in \cD_{<w}(\GrA)$ for some $w \in \bZ$. For example, by Lemma \ref{local_cohom_weight_bounded}, this is true for $M \in \Dperf(\GrA)$ and $N \in \Dbcoh(\GrA)$. It is therefore natural to replace the condition $\RHom_A(M ,\RGam_{I^+}(N)) \simeq 0$ that we considered above by a stronger condition $\RHomcom_A(M ,\RGam_{I^+}(N)) \in \cD_{<0}(\GrA)$. Indeed, this latter condition says that the weight $i = 0$ is already in the ``stable range'' for the pair $(M,N)$, and may be seen as the ``real reason'' why the condition $\RHom_A(M ,\RGam_{I^+}(N)) \simeq 0$ may hold. Accordingly, we make the following

\bdf  \label{strong_section_def}
A full triangulated subcategory $\cE \subset \cD(\GrA)$ is said to be a \emph{strong positive} (resp. \emph{negative}) \emph{slice} if $\RHomcom_A(M ,\RGam_{I^+}(N)) \in \cD_{<0}(\GrA)$ (resp. $\RHomcom_A(M ,\RGam_{I^-}(N)) \in \cD_{>0}(\GrA)$) for all $M,N \in \cE$.
One also speaks of a \emph{strong positive} (resp. \emph{negative}) \emph{slice} of $\cP^{\pm} \subset \cD_{I^{\pm}\text{-triv}}(\GrA)$ as in Definition \ref{partial_sec_def}.
\edf


The technique of weight truncation gives (among other things) a particular recipe for constructing strong positive and negative slices $\cE^+$ and $\cE^-$. In particular, we have the following (see \cite{Yeu20b})
\bpp
The entire category $\cD_{\IpTR}(\GrA)$ admits a (specific) strong positive slice. In particular, there is a fully faithful exact functor $\cL_{[\geq w]} : \cD_{\IpTR}(\GrA) \ra \cD(\GrA)$ such that $\Ce_{I^+} \circ \cL_{[\geq w]} \cong \id$.
Moreover, this functor sends 
$\cD^-_{\coh(\IpTR)}(\GrA)$ to $\Dmcoh(\GrA)$.
The similar statement, with $I^+$ replaced by $I^-$, is true.
\epp

%
%
Roughly speaking, the following analogue of Proposition \ref{RHom_RGam_swap}, with the exact same proof, says that if \eqref{Gor_and_local_cohom_dual_2}(ii) holds, then $\cE^-$ should be at least as large as $\cE^+$.

\bpp 
Assume that the Noetherian $\bZ$-graded ring $A$ satisfies \eqref{Gor_and_local_cohom_dual_2}(ii) for $a \geq 0$. Suppose that $M \in \Dperf(\GrA)$ and $N \in \cD(\GrA)$, then we have 
\begin{equation*}
\RHomcom_A(M,\RGam_{I^+}(N)) \in \cD_{<0}(\GrA) \quad \Rightarrow \quad 
\RHomcom_A(N,\RGam_{I^-}(M)) \in \cD_{>a}(\GrA) \subset \cD_{>0}(\GrA)
\end{equation*}
\epp

\bcor
Denote by $\cD_c \subset \cD$ the subcategory of compact objects.
Suppose that  $\cD_{\IpTR}(\GrA)_c$ has a strong positive slice in $\Dperf(\GrA)$. If $A$ satisfies \eqref{Gor_and_local_cohom_dual_2}(ii) for $a \geq 0$, then there is a fully faithful exact functor
 $\cD_{\IpTR}(\GrA)_c \ra \cD_{\ImTR}(\GrA)_c$
\ecor

\brm  \label{Bondal_Orlov_remark}
A conjecture of Bondal and Orlov says that if $X$ and $X'$ are smooth projective Calabi-Yau varieties that are birational to each other, then they are derived equivalent. 
A result \cite{Kaw08} of Kawamata shows that they are connected by a finite sequence of flops, each of which satisfies the assumptions of Theorem \ref{main_thm_intro}(2) (see Remark \ref{flop_conn_min_rem}). Thus, to prove the Bondal-Orlov conjecture, it suffices to consider the case of a flop in the situation of Theorem \ref{main_thm_intro}(2), where both $X^-$ and $X^+$ are Calabi-Yau, and $X^-$ is smooth. 
In this case, consider the following question:
\begin{equation}  \label{loc_ff_sec_question}
\parbox{40em}{Does the localization $\Ce_{\scI^-} : \Dperf(\GrA) \ra \cD_{\scI^-\text{-triv}}(\GrA)_c$ have a fully faithful section?}
\end{equation}
where we denote by $\cD_c \subset \cD$ the full subcategory of compact objects.

Suppose that \eqref{loc_ff_sec_question} can be answered affirmatively, then take $\cE \subset \Dperf(\GrA)$ to be the essential image of this fully faithful section, so that it is a negative slice in the sense of Definition \ref{partial_sec_def}. By a generalization of Corollary \ref{Dperf_Ce_ff_flop} to the non-affine setting (see Corollary \ref{Dperf_Ce_ff_flop_intro}), it is also a positive slice, and we have a fully faithful functor 
\begin{equation}  \label{Phi_E_ff_functor}
\Phi_{\cE} \, : \, \Dperf(X^-) \simeq \cD_{\scI^-\text{-triv}}(\GrA)_c \xla[\simeq]{\Ce_{\scI^-}} \cE \xrinto{\Ce_{\scI^+}} \cD_{\scI^+\text{-triv}}(\GrA)_c \simeq \Dperf(X^+)
\end{equation}
Now, observe that the canonical DG enhancement of $\Dperf(X^-)$ is smooth; while that of $\Dperf(X^+)$ is proper. Thus, the functor \eqref{Phi_E_ff_functor} has a right adjoint in view of the following result, which follows from \cite[Lemma 2.8]{TV07}:
\begin{equation}  \label{smooth_proper_right_adj}
\parbox{40em}{Let $F : \cA \ra \cB$ be a DG functor between small DG categories, where $\cA$ is smooth and $\cB$ is proper, then the left Kan extension functor $F_! : \Dperf(\cA) \ra \Dperf(\cB)$ has a right adjoint.}
\end{equation}
Since \eqref{Phi_E_ff_functor} is fully faithful, it embeds $\Dperf(X^-)$ as a semi-orthogonal summand of $\Dperf(X^+)$. However, $X^+$ is Gorenstein Calabi-Yau, so $\Dperf(X^+)$ has no non-trivial semi-orthogonal decompositions, and hence $\Phi_{\cE}$ must be an equivalence. Since smoothness of a variety can be detected from its derived category (see, \eg, \cite[Proposition 3.13]{Lun10}), this would show that $X^+$ is smooth as well, and one may proceed to the next flop. Thus, the problem is reduced to answering the formal question \eqref{loc_ff_sec_question}.
The technique of weight truncation gives an affirmative answer to \eqref{loc_ff_sec_question} when $\Spec_Y \cA$ is smooth, or more generally a local complete intersection of a suitable form. It seems possible that \eqref{loc_ff_sec_question} holds quite generally.
\erm

\end{document}